\documentclass[12pt]{amsart}

\usepackage{amssymb}

\usepackage{epsfig,color}

\newtheorem{thm}{Theorem}
\newtheorem{prop}{Proposition}
\newtheorem{defin}{Definition}

\newtheorem{cor}{Corollary}
\newtheorem{lem}{Lemma}
\newtheorem{rem}{Remark}

\newtheorem{exa}{Example}

\newcommand{\N}{{\mathbb N}}

\newcommand{\OO}{{\mathcal O}}

\newcommand{\ppp}{{\mathfrak P}}

\newcommand{\QQ}{{\mathbb Q}}

\newcommand{\R}{{\mathbb R}}\newcommand{\RR}{{\mathbb R}^2}

\newcommand{\Z}{{\mathbb Z}}\newcommand{\ZZ}{\mathbb Z^2}

\newcommand{\sing}{\text{sing}}

\newcommand{\bo}{\partial} %%%%%%%%%% BOUNDARY

\newcommand{\con}{\mbox{cnsv}} %%%%%%%%%% CONSERVATIV
\newcommand{\dis}{\mbox{dspt}}

\newcommand{\gl}{\text{GL}_+(2,\R)}    %%%%%%%%%%%%% GL(2,R)

\newcommand{\hol}{\text{Hol}}\newcommand{\ho}{\text{hol}}

\newcommand{\iso}{\text{Iso}}   %%%%%%%%%%%%%%%% ISOMETRIES
\newcommand{\inter}{\text{interior}}     %%%%%%%%%%%%%%%%% INTERIOR

\newcommand{\mac}{\text{GL}(2,\R)}   %%%%%%%%%%%%% SL(2,R) (FOR McMULLEN)
\newcommand{\vol}{\text{vol}}

\newcommand{\bg}{\bar{g}}

\newcommand{\bet}{\bar{\het}}

\newcommand{\tc}{\tilde{c}}

\newcommand{\tC}{\tilde{C}}

\newcommand{\tG}{\tilde{G}}

\newcommand{\tP}{\tilde{P}}

\newcommand{\tS}{\tilde{S}}
\newcommand{\tT}{\tilde{T}}
\newcommand{\tX}{\tilde{X}}
\newcommand{\tY}{\tilde{Y}}

\newcommand{\tmu}{\tilde{\mu}}
\newcommand{\tnu}{\tilde{\nu}}
\newcommand{\ttau}{\tilde{\tau}}

\newcommand{\tppp}{\tilde{\ppp}}

\newcommand{\al}{\alpha}
\newcommand{\be}{\beta}
\newcommand{\ga}{\gamma}\newcommand{\Ga}{\Gamma}
\newcommand{\de}{\delta}

\newcommand{\ka}{\kappa}
\newcommand{\la}{\lambda}
\newcommand{\Om}{\Omega}
\newcommand{\Si}{\Sigma}
\newcommand{\vp}{\varphi}
\newcommand{\het}{\theta}

\begin{document}

\bibliographystyle{plain}

\title[Noncompact polygonal surfaces]
{Geometry, topology and dynamics of geodesic flows on noncompact
polygonal surfaces}

\author{Eugene Gutkin}

\address{Nicolaus Copernicus University (UMK), Chopina 12/18, Torun 87-100 and
Mathematics Institute of the Polish Academy of Sciences (IMPAN), Sniadeckich 8,
Warszawa 10, Poland} \email{gutkin@mat.umk.pl,gutkin@impan.pl}

%\thanks{The work was partially supported by MNiSzW grant
%NN201384834.}

\keywords{(periodic) polygonal surface, geodesic, skew product,
cross-section, displacement function, recurrence, transience,
ergodicity}

\subjclass{37C40, 37D50, 37E35}

%\date{June 13, 2006}
\date{\today}

\begin{abstract}
We establish the background for the study of geodesics on
noncompact polygonal surfaces. For illustration, we study the
recurrence of geodesics on $\Z$-periodic polygonal surfaces. We
prove, in particular, that almost all geodesics on a topologically
typical $\Z$-periodic surface with a boundary are recurrent.
\end{abstract}

\maketitle

\tableofcontents

\section*{Introduction}       \label{intro}
The billiard ball problem on compact, planar polygons,
$P\subset\RR$, offers a variety of basic open questions. Some of
them have to do with the statistical aspect of dynamics, e. g.,
the ergodicity of polygonal billiards; some others concern the
topological aspect, e. g., the questions about periodic  billiard
orbits in $P$. We refer the reader to \cite{Gut03} for more
detailed information on the subject.

However, the recurrence of billiard orbits on any compact polygon
immediately follows from the poincare recurrence theorem, because
the liouville measure, which is the natural invariant measure,  is
finite. The situation changes, once we pass from compact to
noncompact polygons.  Although the billiard on certain noncompact
polygons has been in the physics literature some hundred years
\cite{Eh},\footnote{This is the English translation; the original
German version appeared in 1912: Encyklop\"{a}die der
Mathematischen Wissenschaften. B.G. Teubner, Leipzig, IV 2, II,
Heft 6. Reprinted in M.J. Klein (ed.): {\em Paul Ehrenfest:
Collected Scientific Papers}, Interscience, New York, 1959. } the
subject has not been systematically studied.

All of the open questions for compact polygons are, of course,
open for noncompact polygons, as well. In addition, because the
liouville measure is infinite, a noncompact polygon,
$P\subset\RR$, may yield a dissipative component in the billiard
dynamics. Thus, we need to decompose the phase space for the
billiard on $P$ into the conservative and dissipative components.

Before we begin the study of the billiard on noncompact polygons,
the question arises: What do we mean by {\em noncompact polygons}?
Indeed, this concept is so large and inhomogeneous that we do not
expect nontrivial results pertaining to the billiard on every
noncompact polygon. In order to obtain them, we should first
identify a class of noncompact polygons which is i) sufficiently
large to deserve a study; ii) sufficiently homogeneous to
(hopefully) possess nontrivial properties that hold for all
polygons in this class. Here we study such a class of noncompact
polygons: The {\em periodic polygons}.

The recurrence, transience and ergodicity of the billiard on
periodic polygons is also studied in \cite{CoGu10}; the present
paper is a supplement to \cite{CoGu10} in a certain sense. We will
now explain in what sense. The study of billiard on noncompact
polygons leads to the concepts of noncompact polygonal surfaces
and their coverings. It also leads to the notion of rationality
for such surfaces, and hence to noncompact translation surfaces.
Their counterparts in the compact case are well known; the basic
properties and relationships between these notions are well
established. They form the background for the study of billiard on
compact polygons. See \cite{Gut86,Gut96} for more on this.

Here we do the same kind of ground work for noncompact polygonal
surfaces. We introduce the basic notions in this subject; in
particular, we introduce and study the coverings of noncompact
polygonal surfaces. We establish basic connections between the
coverings and the holonomy of polygonal surfaces  We define the
rationality for noncompact polygonal surfaces, and introduce
noncompact translation surfaces. This material constitutes the
geometry and topology part of our paper. See section~~\ref{set}.

In section~~\ref{flows} we do the ground work regarding the
geodesic flows on noncompact polygonal surfaces. Let $P$ be one.
If $P$ is rational, the geodesic flow on $P$ decomposes into a
one-parameter family of directional flows. These can be realized
as the linear flows on the noncompact translation surface
$S=S(P)$. In the subject of compact, rational polygonal surfaces,
{\em arithmetic surfaces} play a special role \cite{Gut84,GJ96}.
In section~~\ref{arithm_sub} we introduce and study noncompact
arithmetic polygonal surfaces.

The material of section~~\ref{set} and section~~\ref{flows} is
used in \cite{CoGu10} to study the recurrence and ergodicity for
geodesic dynamics on periodic surfaces. In particular,
\cite{CoGu10} contains results on ergodicity for certain
$\ZZ$-periodic surfaces. Here, we restrict ourselves to the
simpler class of $\Z$-periodic  polygonal surfaces. We study in
detail the recurrence of geodesics on these surfaces.
Concentrating on $\Z$-periodic surfaces with boundaries, we
establish a few results on recurrence of the geodesic flow and the
directional flows. See Theorem~~\ref{per_rat_recurr_thm} and
Theorems~~\ref{dens_gdel_recur_thm},~~\ref{mor_dens_gdel_thm}.
Section~~\ref{transi} represents the dynamics part of our paper.

\medskip

As we already mentioned, the goal in the present study is to
establish the background for the subject of noncompact polygonal
surfaces. This determined the style of the paper, which differs
somewhat from that of a typical research paper. Our emphasis is on
basic principles and the motivations, rather than on elaborate
mathematical results. Thus, we present relatively few theorems,
and offer many examples. The examples serve to illustrate the
basic notions and to motivate the definitions.

\medskip

\noindent{\bf Acknowledgements.} The present paper owes much to the 
collaboration \cite{CoGu10} with Jean-Pierre Conze. The author 
gave several presentations of this work, in particular, at
the ICTP Conference on Pseudo-Chaos in Trieste in September 
2009 and at the Seminar on Mathematical Physics at UCI in February 2010.
The audience's feedback is greatfully acknowledged. Especially useful
were the comments of Carl Dettmann and Lana Jitomirskaya.
The work was partially supported by MNiSzW grant
NN201384834.

\section{The setting and preliminaries}       \label{set}
In this section we establish the subject of our study: Noncompact
polygonal surfaces. For readers' convenience, we will first recall
the basic material on compact polygonal surfaces.
\subsection{Compact polygonal surfaces}   \label{comp_sub}
\hfill \break We begin with the basic example of a compact
polygonal surface: A planar polygon. By this we mean a compact
region, $P\subset\RR$, whose boundary consists of a finite number
of linear intervals. We view $P$ as a riemannian manifold with
boundary and corners. Then the geodesics are the billiard orbits
in $P$; the geodesic flow is the billiard flow. Billiard orbits in
polygons are of great interest to physicists \cite{Gut86,Gut96};
they have done extensive numerical studies of polygonal billiards.
Many mathematical results in this subject concern {\em rational
polygons} \cite{Smi00}. A polygon is rational if the angles
between its sides are rational multiples of $\pi$. Denote by
$\hol_r(P)\subset O(2)$ the linear part of the group generated by
the orthogonal reflections with respect to the sides of
$P$.\footnote{The meaning of this notation will become clear
later on.} Then $P$ is rational iff $\hol_r(P)=D_N$, the dihedral
group of order $2N$, where $N=N(P)$ is determined by the angles of
$P$. For irrational polygons $|\hol_r(P)|=\infty$.

By a surface we will always mean a connected, orientable surface.
This requirement is not essential: We impose it to simplify the
exposition. Any compact polygonal surface, say $S$, can be
decomposed into a finite number of polygons, say $P_1\cup\dots\cup
P_t$. To recover $S$ from $P_1,\dots, P_t$, we glue them along
(some of) their sides by isometries. This is the {\em defining
mapping}, $f:\cup_{1\le i \le t}P_i\to S$, of the disjoint union
of comprising polygons onto the surface. As a space, $S$ is a
two-dimensional manifold, in general, with a boundary. The
boundary, $\bo S$ is the union of the sides of $P_1,\dots, P_t$
that are left unglued. The euclidean structure of $P_1,\dots, P_t$
endows $S$ with a locally euclidean metric. The metric may be
singular only at the points of $S$ which come from the corners of
$P_1,\dots, P_t$. Thus, the singular set, $\Si\subset S$, is
finite. Near a point $c\in\Si\cap\inter(S)$ (resp. $c\in\Si\cap\bo
S$) the surface is isometric to a euclidean cone (resp. euclidean
wedge). The apex angle of this cone (resp. wedge) is the {\em cone
angle} (resp. {\em wedge angle}) of the {\em cone point} (resp.
{\em wedge point}). Cone (resp. wedge) angles of singular
points may take any values except $2\pi$ (resp. $\pi$). Compact
polygonal surfaces without singular points are either flat tori or
flat cylinders.

\begin{exa}      \label{pillow_exa}
{\em Let $P',P''$ be two copies of a simple $k$-gon $P=A_1\dots
A_k$; let $\al_1,\dots,\al_k$ be the respective angles. Let $S$ be
the polygonal surface formed by gluing the corresponding sides
of $P'$ and $P''$ reversing the orientation. This surface is the
{\em doubling of $P$}. We will denote it by $DP$. Topologically,
$DP\sim S^2$. Its singular set consists of $k$ cone points, with
cone angles $2\al_1,\dots,2\al_k$.

}
\end{exa}

Let $\iso(\RR)$ (resp. $\iso_0(\RR)$) be the group of (resp.
orientation preserving) isometries. We recall the notions of
$(G,X)$-manifolds and $(G,X)$-orbifolds of Ehresmann and Thurston
\cite{Th97}. Heuristically, compact polygonal surfaces are
$(\iso(\RR),\RR)$-manifolds with singularities. In particular, we
associate with a compact polygonal surface $S$ its {\em holonomy
group} $\hol(S)\subset\iso(\RR)$. We define $\hol(S)$ the same way
it is defined for $(\iso(\RR),\RR)$-orbifolds. Namely, we choose a
base point $s_0\in S\setminus\Si$, and consider {\em regular} (i.
e., avoiding the singularities) {\em loops}  based at $s_0$. If
$\ga$ is such a loop, we {\em develop} $S$ along $\ga$ yielding
$h(\ga)\in\iso(\RR)$ which depends only on the regular homotopy
class of $\ga$. The preceding construction is the {\em holonomy
homomorphism} $\ho:\pi_1(S\setminus\Si)\to\iso(\RR)$; its range is
the holonomy group $\hol(S)$.

The semi-direct product decomposition $\iso(\RR)=O(2)\times\RR$
yields $\hol(S)=\hol_r(S)\times\hol_t(S)$, where $\hol_r(S)$
(resp. $\hol_t(S)$) is the  {\em rotational (resp. translational)
holonomy}. We say that $S$ is a {\em rational polygonal surface}
if $|\hol_r(S)|<\infty$. If $P$ is a polygon, $\hol_r(P)$ is
generated by linear orthogonal reflections about its sides. Thus,
for polygons the present definition of ``rational'' agrees with
the earlier one. We have
$\hol(DP)=\hol(P)\cap\iso_0(\RR)$.\footnote{This is a special case
of a general statement about the holonomy of coverings, see
below.} Hence, the doubling surface $DP$ in
Example~~\ref{pillow_exa} is rational iff $P$ is a rational
polygon.

\begin{defin}      \label{covering_def}
{\rm Let $R,S$ be compact polygonal surfaces. A {\em covering of polygonal surfaces}
is a continuous, surjective mapping $\vp:R\to S$ compatible with the
respective $(\iso(\RR),\RR)$-structures.

}
\end{defin}
\begin{exa}      \label{more_pillow_exa}
{\em Let $S$ be any compact polygonal surface with a boundary.
Gluing up the two copies $S',S''$ along the boundaries, we obtain
a closed polygonal surface $DS$, the {\em doubling of $S$}. When
$S$ is a polygon, we recover Example~~\ref{pillow_exa}. The
natural projection $\vp:DS\to S$ is a covering of polygonal
surfaces.

}
\end{exa}
\begin{rem}    \label{double_rem}
{\em The trick of doubling allows us to restrict our
considerations, if need be, to boundaryless  polygonal surfaces.

}
\end{rem}
\begin{prop}  \label{cover_hol_prop}
1. Let $\vp:R\to S$ be a covering of compact polygonal surfaces.
It induces an inclusion $\hol(R)\subset\hol(S)$, compatible with
the injection
$\vp_*:\pi_1(R\setminus\Si(R))\to\pi_1(S\setminus\Si(S))$ and with
the holonomy homomorphisms
$\ho_R:\pi_1(R\setminus\Si(R))\to\iso(\RR),
\ho_S:\pi_1(S\setminus\Si(S))\to\iso(\RR)$.

2. Let $S$ be a compact polygonal surface. Let $H\subset\hol(S)$
be a subgroup of finite index. Then there is a unique covering of
compact polygonal surfaces $\vp_H:\tS_H\to S$ such that
$\hol(\tS_H)=H$.
\begin{proof}
Choose regular base points $r_0\in R,s_0\in S$ so  that
$\vp(r_0)=s_0$. Let $\al$ be a regular loop in $(R,r_0)$. Then
$\be=\vp(\al)$ is a regular loop in $(S,s_0)$. We denote by
$[\al],[\be]$ their homotopy classes. Simultaneously developing
$R$ along $\al$ and $S$ along $\be$, we obtain
$\ho_R(\al)=\ho_S(\be)$. Replacing loops by their homotopy classes
and using that $[\be]=\vp_*([\al])$, we obtain the first claim.

For the proof of the second claim we make use of
Remark~~\ref{double_rem} and assume that $S$ is a closed surface.
Set $G=\ho_S^{-1}(H)\subset\pi_1(S\setminus\Si(S))$. The index
$[G:\pi_1(S\setminus\Si(S))]$ is equal to the index of $H$ in
$\hol(S)$. By assumption, $[G:\pi_1(S\setminus\Si(S))]=d\in\N$.
Let $\vp_0:R_0\to S\setminus\Si(S)$ be the corresponding
topological covering of degree $d$. The
$(\iso(\RR),\RR)$-structure pulls back from $S\setminus\Si(S)$ to
$R_0$. The covering $\vp_0:R_0\to S\setminus\Si(S)$ uniquely
extends to a projection $\vp:R\to S$ of their completions; the
mapping $\vp:R\to S$ is a branched covering of degree $d$. The
branching locus is $\Si(S)$. Thus, $R$ is a compact polygonal
surface.
Set $\tS_H=R$. By construction, $\vp_H:\tS_H\to S$ is a covering
of polygonal surfaces. It is straightforward to check that it has
the required properties.
\end{proof}
\end{prop}

\begin{cor}   \label{ratio_cov_cor}
Let $\vp:R\to S$ be a covering of compact polygonal surfaces. Then
one of the surfaces is rational iff both are rational.
\begin{proof}
By Proposition~~\ref{cover_hol_prop}, $\hol_r(R)\subset\hol_r(S)$.
Thus, if $S$ is a rational surface, then $|\hol_r(R)|<\infty$, i.
e., $R$ is a rational surface. Suppose now that $S$ is irrational,
i. e., $|\hol_r(S)|=\infty$. By Proposition~~\ref{cover_hol_prop},
$[\hol_r(R):\hol_r(S)]<\infty$, implying $|\hol_r(R)|=\infty$.
\end{proof}
\end{cor}

A {\em compact translation surface} is a compact polygonal surface
whose rotational holonomy  is trivial. Equivalently, a compact
translation surface is a compact polygonal surface which carries a
$(\RR,\RR)$-structure \cite{GJ96,GJ00}. Also equivalently, a
compact polygonal surface $S$ is a translation surface iff
$\hol(S)=\hol_t(S)$. Compact translation surfaces arise in several
contexts, e. g., in complex analysis. They are instrumental in the
analysis of billiard in rational polygons \cite{Smi00}. A {\em
covering of translation surfaces} is a covering of polygonal
surfaces $\vp:R\to S$, where $R$ and $S$ are translation surfaces.

\begin{cor}   \label{trans_cov_cor}
Let $S$ be a compact, rational polygonal surface. Then it has a
unique minimal covering $\vp_t:R\to S$ by a compact translation
surface. The minimality means that if $\vp:R'\to S$ is any
covering by a compact translation surface then there is a covering
of translation surfaces $\psi:R'\to R$ such that
$\vp=\vp_t\circ\psi$.
\begin{proof}
We have $[\hol_t(S):\hol(S)]=|\hol_r(S)|<\infty$. Set
$H=\hol_t(S)$; let $R=\tS_H$ and let $\vp_t:R\to S$ be the covering
in Proposition~~\ref{cover_hol_prop}. Then
$\hol_t(R)=\hol(R)=\hol_t(S)$, i. e., $R$ is a translation
surface. We will now prove the minimality of $\vp_t:R\to S$. If
$\vp:R'\to S$ is any covering by a translation surface, then
$\hol(R')\subset\hol_t(S)=\hol(R)$. Let $\psi:R'\to R$ be the
covering constructed in Proposition~~\ref{cover_hol_prop}. Then
$\psi:R'\to R$ is a covering of polygonal surfaces. By
construction, $\vp=\vp_t\circ\psi$.
\end{proof}
\end{cor}

Let $P$ be any compact, rational polygonal surface. Let
$\vp_t:S(P)\to P$ be the minimal covering ensured by
Corollary~~\ref{trans_cov_cor}. We will refer to $S(P)$ as the
(canonical) translation surface of $P$ and to $\vp_t:S(P)\to P$ as
the {\em canonical translation covering}. Note that the degree of
$\vp_t:S(P)\to P$ is $|\hol_r(P)|$. If $P$ is a rational polygon,
$S(P)$ is often called the {\em Katok-Zemlyakov surface}. The term
acknowledges the work \cite{KZ} which derived nontrivial dynamical
consequences from $\vp_t:S(P)\to P$. We point out that these
coverings, for any compact rational polygonal surface, are in the
literature since 1907. See the references in \cite{Gut84}.

\medskip

\subsection{Noncompact polygons: Definitions and examples}   \label{noncompol_sub}
\hfill \break Our goal is to extend the above material to
noncompact polygonal surfaces. As in section~~\ref{comp_sub}, the
basic example of such surface is a {\em noncompact polygon}. By
this we mean a closed, noncompact region, $\Om\subset\RR$ such
that $\bo\Om$ consists of (maximal) line segments: The sides of
$\Om$. They may have finite or infinite length; there may be
infinitely many sides. We will introduce some requirements on our
noncompact polygons.

\begin{itemize}

\item[ A.] Any set $\{|x|,|y|\le a\}$ intersects at most a finite number of sides.

\item[ B.] The side lengths are bounded away from zero.

\end{itemize}

\begin{rem}      \label{no_snow_rem}
{\em In view of condition A, we exclude from consideration {\em
fractal polygons}. We feel that the billiard on fractal polygons
and {\em fractal polygonal surfaces} is a separate subject.
Condition B is not as paramount as condition A; some natural
noncompact polygons do not satisfy it. See
Example~~\ref{half_stair_exa} below. We point out that
regions,  such as the whole plane, half-planes, wedges, infinite bands, etc are
noncompact polygons.

}
\end{rem}

We will now discuss a few examples.

\begin{exa}      \label{fin_obst_exa}
{\em Let $P_1,P_2,\dots,P_t\subset\RR$ be simple polygons, whose
interiors are pairwise disjoint. Set
$\Om=\RR\setminus\inter(P_1\cup P_2\cup\cdots\cup P_t)$. Then
$\Om$ is a noncompact polygon. We will refer to it as the plane
with (a finite number of) {\em polygonal obstacles}.

}
\end{exa}

\begin{figure}[htbp]
\begin{center}
\input{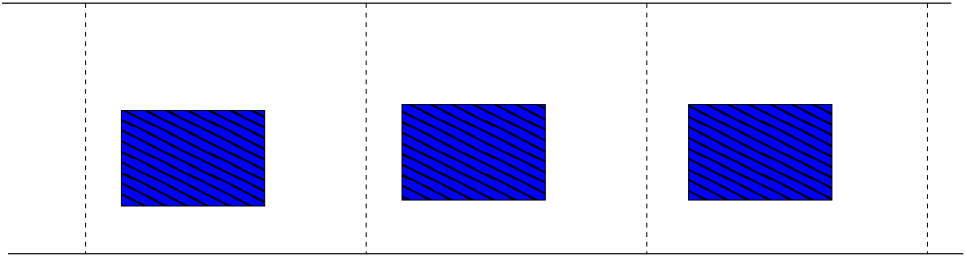}
\caption{An infinite band with a $\Z$-periodic configuration of
rectangular obstacles.} \label{rect_loren_band}
\end{center}
\end{figure}

\begin{exa}      \label{period_obst_exa}
{\em By an {\em infinite band} we will mean the planar region
bounded by two parallel lines. The {\em standard band} $B_0$
is bounded by $\{y=0\}$ and $\{y=1\}$. By a {\em standard $a\times
b$ rectangle} we mean any rectangle
$R(a,b;\xi,\eta)=\{(x,y):\xi\le x\le\xi+a,\eta\le y\le\eta+b\}$.
Setting $\xi=\eta=0,a=b=1$, we obtain the {\em standard unit
rectangle} $R$. We set $R_{(m,n)}=R+(m,n)$.

Let $P\subset\inter(R)$ be a simple polygon. The region
$\Om=B_0\setminus\{\cup_{k\in\Z}(P+(k,0))\}$ is a noncompact
polygon. This is an infinite {\em band with a $\Z$-periodic configuration of
polygonal obstacles}. Figure~~\ref{rect_loren_band} shows a
special case.

}
\end{exa}

\begin{figure}[htbp]
\begin{center}
\input{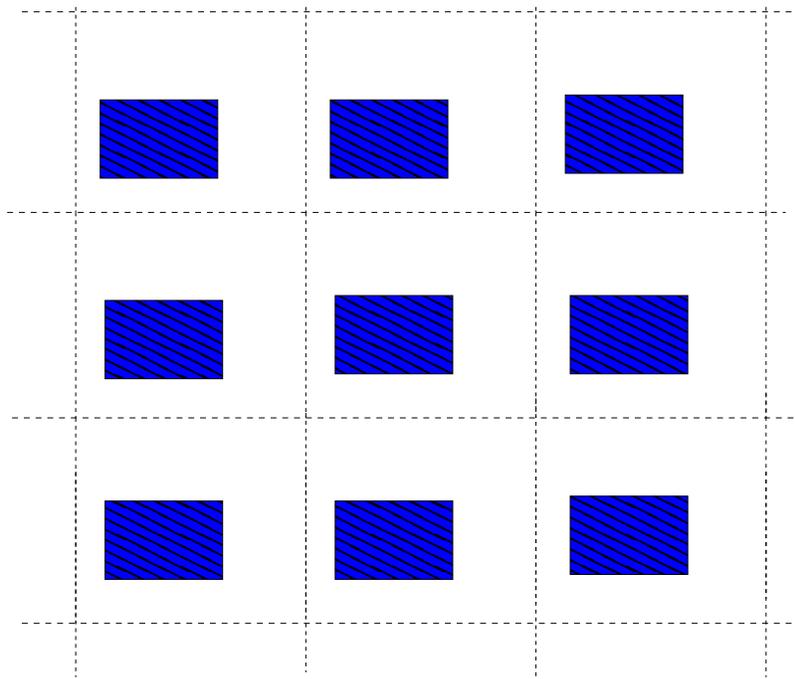}
\caption{The euclidean plane with a $\ZZ$-periodic configuration
of rectangular obstacles.} \label{rect_loren_gas}
\end{center}
\end{figure}

\begin{exa}      \label{doubl_period_obst_exa}
{\em Let $P$ be a simple polygon, as in
Example~~\ref{period_obst_exa}. Set
$\Om=\RR\setminus\{\cup_{(m,n)\in\ZZ}(P+(m,n))\}$. This noncompact
polygon is the {\em plane with a $\ZZ$-periodic configuration of
polygonal obstacles}. The special case, when $P$ is a
rectangle, has been in the literature for some hundred years
\cite{Eh}. Following
\cite{Eh} and \cite{HW80}, this polygonal surface is often called
the {\em wind-tree model}. Figure~~\ref{rect_loren_gas} shows a
$\ZZ$-periodic configuration of rectangular obstacles.

}
\end{exa}

\medskip

\begin{figure}[htbp]
\begin{center}
\input{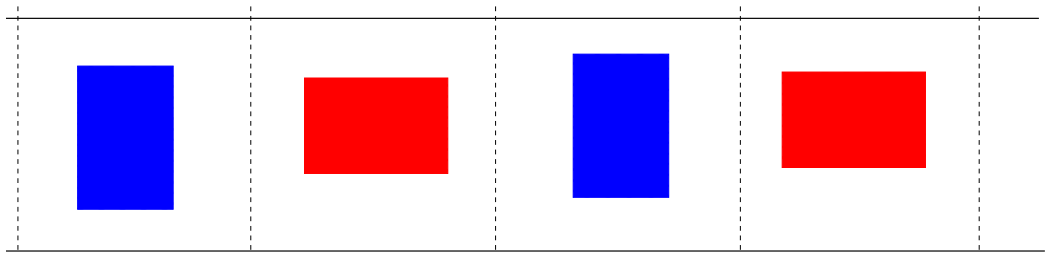}
\caption{A band with a $\Z$-quasi-periodic (or $\Z$-periodic)
configuration of rectangular obstacles.} \label{quas_rect_band}
\end{center}
\end{figure}

\begin{figure}[htbp]
\begin{center}
\input{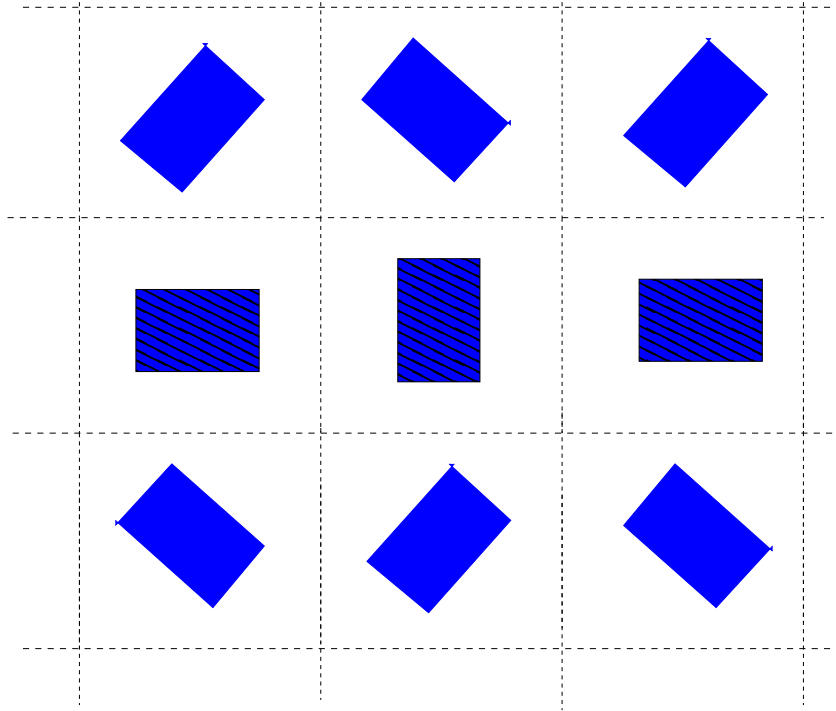}
\caption{The euclidean plane with a $\ZZ$-quasi-periodic (or
$\ZZ$-periodic) configuration of rectangular obstacles.}
\label{quas_rect_gas}
\end{center}
\end{figure}

\begin{exa}      \label{qua_perio_obst_exa}
{\em We will denote by $\rho_{\het}$ the rotation of $\RR$ by the
angle $\het$ around any center point. Let
$P\subset\inter(R)$ be a simple polygon with a
distinguished point, $o\in P$, which will be our center point. We
assume that $\rho_{\het}(P)\subset\inter(R)$ for any
$\het$.

Let $0<\al,\be<2\pi$. The region
$\Om=B_0\setminus\{\cup_{k\in\Z}(\rho_{k\al}(P)+(k,0))\}$ is a
noncompact polygon. If $\al/\pi$ is rational (resp. irrational),
$\Om$ is the standard band with a $\Z$-periodic (resp. {\em
$\Z$-quasi-periodic}) configuration of polygonal obstacles. We
will say that $\Om$ is a {\em periodic or a quasi-periodic band},
for brevity. Figure~~\ref{quas_rect_band} shows an example.

Let
$\Om=\RR\setminus\{\cup_{(m,n)\in\ZZ}(\rho_{m\al+n\be}(P)+(m,n))\}$.
This noncompact polygon is the {\em plane with a $\ZZ$-periodic or
$\ZZ$-quasi-periodic configuration of polygonal obstacles},
depending on $\al,\be$. See figure~~\ref{quas_rect_gas} for an
example.

}
\end{exa}

\medskip

Preceding examples satisfied conditions A and B. Our next example
is a class of noncompact polygons satisfying A but not necessarily
B.

\begin{figure}[htbp]
\begin{center}
\input{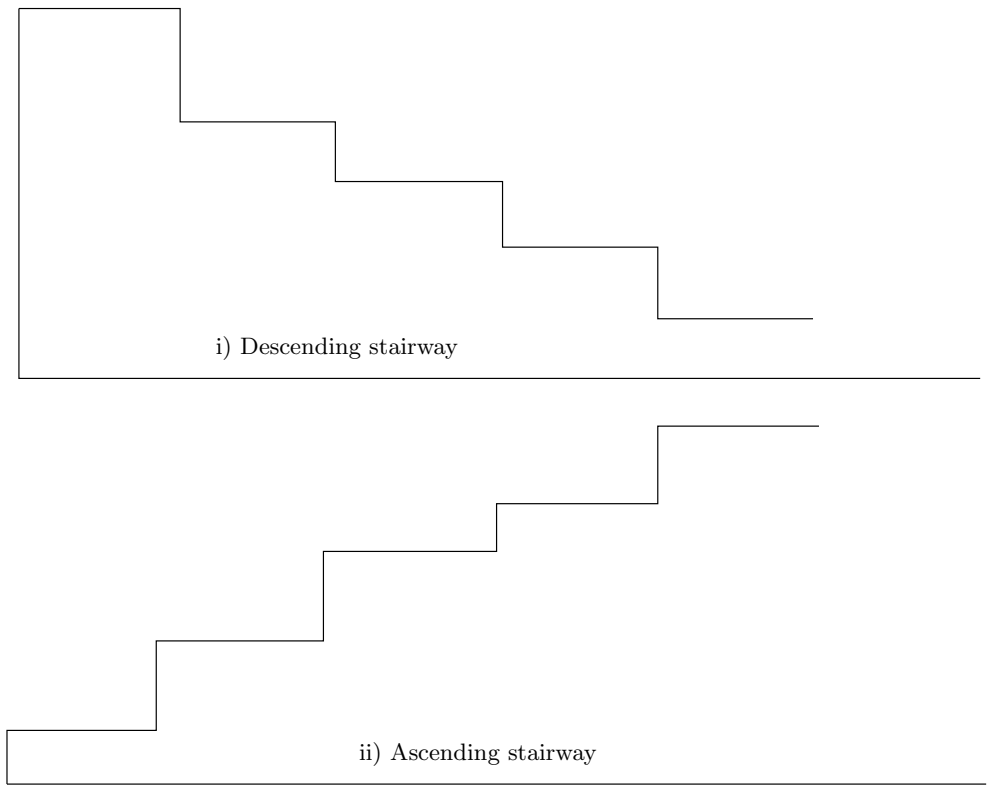}
\caption{Examples of stairways: i) descending, ii) ascending.}
\label{stairways}
\end{center}
\end{figure}

\begin{figure}[htbp]
\begin{center}
\input{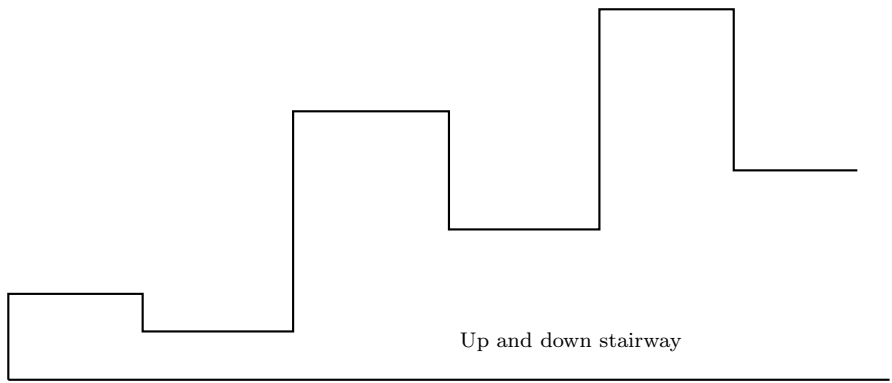}
\caption{An up and down stairway.} \label{upanddown}
\end{center}
\end{figure}

\begin{exa}      \label{half_stair_exa}
{\em Let $(h)=h_0,h_1,\dots$ be an infinite sequence of (strictly)
positive numbers. The noncompact polygons $P(h)=\cup_{k\ge
0}R(1,h_k;k,0)$ are  {\em infinite stairway polygons}, or {\em
stairways} for brevity. The work \cite{Mirko} studied the billiard
on $P(h)$ when the sequence $(h)=h_0,h_1,\dots$ is strictly
decreasing. We call the polygons $P(h)$ with $(h)$ strictly
decreasing (resp. increasing) the {\em descending (resp.
ascending) stairways}. If the sequence $(h)$ is not monotone, we
will informally say that $P(h)$ is an {\em up and down stairway}.
See figures~~\ref{stairways} and~~\ref{upanddown}  for an
illustration.

}
\end{exa}

If $\Om$ is a noncompact polygon, we denote by
$\hol(\Om)\subset\iso(\RR)$ the group generated by orthogonal
reflections about the lines extending the sides of $\Om$; it is a
semi-direct product $\hol(\Om)=\hol_r(\Om)\times\hol_t(\Om)$. The
polygon $\Om$ is rational if $|\hol_r(\Om)|<\infty$. Then
$\hol_r(\Om)=D_N$; compare with section~~\ref{comp_sub}. If $\Om$
is irrational, the group $\hol_r(\Om)$ may be infinitely
generated.

In Example~~\ref{fin_obst_exa}, $\Om$ is rational iff for any pair
$1\le i,j\le t$ the angles between the sides of $P_i$ and $P_j$
are $\pi$-rational. In Example~~\ref{period_obst_exa}, $\Om$ is
rational iff the angles of $P$ and the angles between the sides of
$P$ and the horizontal axis are  $\pi$-rational. In particular,
the noncompact polygons in figures~~\ref{rect_loren_band}
and~~\ref{rect_loren_gas} are rational, and $N=2$ for both. The
infinite band with rectangular obstacles in
figure~~\ref{quas_rect_band} is rational iff $\al\in\pi\QQ$, iff
$\Om$ is $\Z$-periodic. The polygon $\Om$ in
figure~~\ref{quas_rect_gas} is rational iff $\al,\be\in\pi\QQ$,
iff $\Om$ is $\ZZ$-periodic. The stairway polygons in
Example~~\ref{half_stair_exa} are rational, and $N=2$.

%A noncompact polygon is {\em tame} if the angles between its
%non-parallel sides are bounded away from zero. This is a
%particular case of a concept that plays an important role for
%noncompact polygonal surfaces. See Definition~~\ref{tame_def}
%below.

%\medskip

%
\subsection{Noncompact polygonal surfaces: definitions and examples}   \label{noncomsurf_sub}
\hfill \break A noncompact polygonal surface is assembled from a
disjoint union of a countable (at most) collection of compact and
noncompact polygons. All boundary identifications are via
isometries on full sides of these polygons. Let $S$ be a polygonal surface,
let $\cup_{i\in I}P_i$ be the disjoint union of polygons
comprising $S$, and let $f:\cup_{i\in I}P_i\to S$ be the {\em
defining mapping}.  We impose the following restriction.

\begin{itemize}

\item[ C.] For any $c\in S$ we have $|f^{-1}(c)|<\infty$.

\end{itemize}

The boundary $\bo S$ is the union of the sides of polygons
$P_i,i\in I,$ that are left unglued. The euclidean structure of
polygons $P_i,i\in I,$ endows $S$ with a locally euclidean metric; it may
be singular only at the points which come from the corners of
polygons $P_i,i\in I$. Let $\Si'\subset S$ be the $f$-range of the
set of corners. Since, by condition A, the set $\Si'$ is discrete,
the singular set $\Si\subset\Si'$ is discrete, as well. In view of
condition C, near a point $c\in\Si'\cap\inter(S)$ (resp.
$c\in\Si'\cap\bo S$) the surface is isometric to a euclidean cone
(resp. euclidean wedge). The singular set $\Si$ consists of those
points $c\in\Si'\cap\inter(S),c\in\Si'\cap\bo S$ whose {\em cone
angle} (resp. {\em wedge angle}) is not $2\pi$ (resp. not $\pi$). If
condition B is not satisfied, then $\Si$ may contain infinite
sequences $a_k\ne b_k$ such that the distance $d(a_k,b_k)$ goes to
zero. We finish this section with a few examples.

%By the preceding discussion, a noncompact polygonal surface is a
%locally euclidean surface with a boundary, cone points, and
%corners. Since $S\setminus\Si$ is a $(\iso(\RR),\RR)$-manifold, it
%has a {\em holonomy group} $\hol(S)\subset\iso(\RR)$. We will
%elaborate on this in section~~\ref{noncom_holon_sub}.

%
\begin{exa}      \label{inf_pillow_exa}
{\em Let $\Om$ be a noncompact polygon. Its {\em doubling} $D\Om$
is the noncompact polygonal surface without boundary obtained by
identifying the respective sides of two copies of $\Om$. This is
the noncompact version of Example~~\ref{pillow_exa}.

If $\Om$ is the half-plane, then $D\Om=\RR$. If $\Om$ is the wedge
with wedge angle $\al$, then $D\Om$ is the euclidean cone with the
cone angle $2\al$. Let $\Om=\RR\setminus\inter(P_1\cup
P_2\cup\cdots\cup P_t)$ be as in Example~~\ref{fin_obst_exa}. The
topology of the surface $D\Om$ depends only on $t$. If $t=1$ then
$D\Om\sim\R\times S^1$ is a topological cylinder. The geometry of
$D\Om$ does depend on the polygons $P_1,P_2,\dots,P_t$. Let
$t=1,P_1=P$, and let $\al_1,\dots,\al_p$ be the interior angles of
$P$. Then $D\Om$ has $p$ cone points, with cone angles
$4\pi-2\al_1,\dots,4\pi-2\al_p$ respectively. We leave the
analysis of $D\Om$ for $t>1$ to the reader.

}
\end{exa}

Our next example is a variation on the theme of
Example~~\ref{inf_pillow_exa}.

\begin{figure}[htbp]
\begin{center}
\input{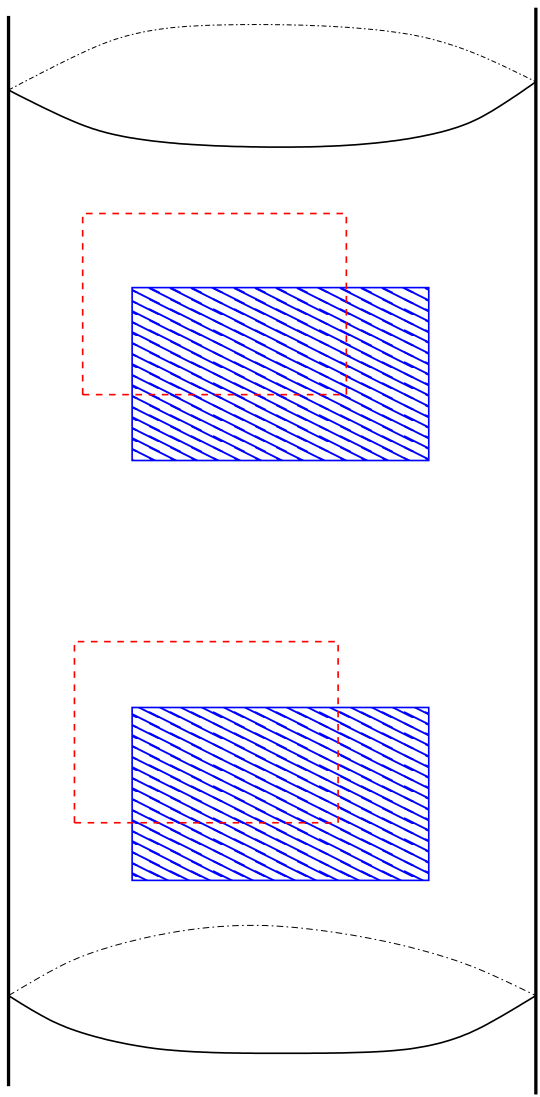}
\caption{Polygonal surface in the shape of a tower with a
$\Z$-periodic collection of pairs of rectangular windows.}
\label{inf_cyl_wind}
\end{center}
\end{figure}

\begin{exa}      \label{hole_pil_exa}
{\em Let $\Om$ be a noncompact polygon, and let $L$ be the
collection of its sides, i. e., $\bo\Om=\cup_{l\in L}l$. Let
$L=G\cup O$ be a partition, with $G\ne\emptyset$. Let $\Om',\Om''$
be two copies of $\Om$, and let $\bo\Om',\bo\Om''$ be the
respective boundaries. For $l\in L$ let $l',l''$ be the
corresponding pair of sides. Let $D_{O}\Om$ be the union
$\Om'\cup\Om''$, where $l'$ and $l''$ are identified iff $l\in G$.
Then $D_{O}\Om$ is a noncompact, connected polygonal surface. We
have $\bo(D_{O}\Om)=\cup_{l\in O}(l'\cup l'')$. In particular,
$D_{\emptyset}\Om=D\Om$ is the only case when
$\bo(D_{O}\Om)=\emptyset$.

Let $\Om$ be the band in figure~~\ref{rect_loren_band}. Let $G$
(resp. $O$) be the pair of infinite segments (resp. the collection
of finite segments)  in $\bo\Om$. Then $D_{O}\Om$ is the infinite
flat cylinder with $\Z$-periodic collections of rectangular
obstacles on its front and its back. This surface resembles a
building in the shape of a tower with infinitely many floors; each
floor of this tower has two identical rectangular windows, one in
the front and the other in the back. See
figure~~\ref{inf_cyl_wind}.  Let $P_0,P_1,\dots,P_t$ be simple,
disjoint polygons. Let $\Om=\RR\setminus\inter(P_0\cup
P_1\cup\dots\cup P_t)$ be the polygon in
Example~~\ref{fin_obst_exa}. Let $G$ be the collection of sides of
$P_0$. Thus, $O$ consists of the sides of $P_1,\dots,P_t$. The
polygonal surface $D_O\Om$ is the flat cylinder $D(\RR\setminus
P_0)$ with the polygonal obstacles $P_1'\cup
P_1''\cup\cdots\cup P_t'\cup P_t''$.

}
\end{exa}

\subsection{Holonomy and coverings of noncompact polygonal surfaces}   \label{noncom_holon_sub}
%
%\hfill \break
Definitions of the holonomy group $\hol(S)\subset\iso(\RR)$ and
the holonomy homomorphism $\ho:\pi_1(S\setminus\Si)\to\iso(\RR)$
of section~~\ref{comp_sub} directly extend to noncompact polygonal
surfaces. As in section~~\ref{comp_sub}, we have a semidirect
product $\hol(S)=\hol_r(S)\times\hol_t(S)$, where the group
$\hol_r(S)\subset O(2)$ (resp. $\hol_t(S)\subset\RR$) is the
rotational (resp. translational) holonomy. These groups may be
infinitely generated. Thus, let $P$ be a noncompact polygon with
infinitely many sides. The group $\hol(P)$ (resp. $\hol_r(P)$) is
generated by (resp. linear parts of) orthogonal reflections about
the lines extending the sides of $P$. Obvious modifications of
Example~~\ref{period_obst_exa} and
Example~~\ref{doubl_period_obst_exa} yield polygons with
infinitely generated groups $\hol_r(P)$, and hence $\hol(P)$. Since
$\hol(DP)=\hol(P)\cap\iso_0(\RR)$, their doublings are examples
of boundaryless noncompact polygonal surfaces with infinitely
generated holonomy groups.

%A cone (resp. wedge) point $c\in\Si(S)$ with cone (resp. wedge)
%angle $\al$ (resp. $\be$) yields the subgroup in $\hol_r(S)$ of
%rotations by $k\al,k\in\Z$ (resp. $2k\be,k\in\Z$).
The concept of {\em covering of polygonal surfaces}, as stated in
Definition~~\ref{covering_def} for compact surfaces, applies to
arbitrary polygonal surfaces. The following lemma summarizes the
basic properties of coverings. See
Proposition~~\ref{cover_hol_prop} for a proof.

\begin{lem}  \label{noncom_cover_lem}
Let $\vp:R\to S$ be a covering of polygonal surfaces. It induces
an inclusion $\hol(R)\subset\hol(S)$, compatible with the
injection
$\vp_*:\pi_1(R\setminus\Si(R))\to\pi_1(S\setminus\Si(S))$ and with
the holonomy homomorphisms
$\ho_R:\pi_1(R\setminus\Si(R))\to\iso(\RR),
\ho_S:\pi_1(S\setminus\Si(S))\to\iso(\RR)$.
\end{lem}
\begin{defin}      \label{tame_cov_def}
{\rm A covering of polygonal surfaces $\vp:R\to S$ is {\em tame}
if the index $[\hol_r(R):\hol_r(S)]<\infty$.

}
\end{defin}
\begin{prop}  \label{noncom_cover_prop}
Let $S$ be a polygonal surface. Let $H\subset\hol(S)$ be a
subgroup. Let $H=H_r\times H_t$ be the decomposition corresponding
to $\hol(S)=\hol_r(S)\times\hol_t(S)$.

Let $[H_r:\hol_r(S)]<\infty$. Then there is a unique covering of
polygonal surfaces $\vp_H:\tS_H\to S$ such that $\hol(\tS_H)=H$.
Its degree is equal to the index of $H$ in $\hol(S)$.
\begin{proof}
In view of the doubling construction of
Example~~\ref{inf_pillow_exa}, we restrict the discussion to
surfaces without boundary. Set
$G=\ho_S^{-1}(H)\subset\pi_1(S\setminus\Si(S))$. Let $\vp_0:R_0\to
S\setminus\Si(S)$ be the unique topological covering corresponding
to $G$. Pulling back by $\vp_0$ the $(\iso(\RR),\RR)$-structure on
$S\setminus\Si(S)$, we endow $R_0$ with a
$(\iso(\RR),\RR)$-structure.

Let $R$ be the completion of $R_0$ with respect to the induced
metric. The assumption $[H_r:\hol_r(S)]<\infty$ ensures that
$R\setminus R_0$ consists of cone points. The covering
$\vp_0:R_0\to S\setminus\Si(S)$ uniquely extends to a branched
covering $\vp:R\to S$ whose ramification locus is contained in
$R\setminus R_0$. The ramification number $r(\tc)$ at a point
$\tc\in R\setminus R_0$ is less than or equal to the index
$[H_r:\hol_r(S)]$. Let $c=\vp(\tc)$ and let $\al(c),\al(\tc)$ be
the respective cone angles. Then
\begin{equation}   \label{bound_eq}
\al(c)\le\al(\tc)=r(\tc)\al(c)\le[H_r:\hol_r(S)]\al(c).
\end{equation}
Thus, $R$ is a polygonal surface. Set $R=\tS_H,\vp=\vp_H$. The
equality $\hol(\tS_H)=H$ holds by construction. The reader will
easily check the remaining claims.
\end{proof}
\end{prop}
\begin{defin}      \label{rational_def}
{\rm A polygonal surface $S$ is {\em rational} (resp. {\em
translation surface}) if $|\hol_r(S)|<\infty$ (resp.
$|\hol_r(S)|=1$).

}
\end{defin}

A polygon is {\em rational} if it is a rational polygonal surface.
Note that a compact polygon is rational in the sense of
Definition~~\ref{rational_def} iff it is rational in the sense of
section~~\ref{comp_sub}. A noncompact polygon, say $P\subset\RR$,
is rational iff there exists $N\in\N$ such that all angles between
the sides of $P$ belong to the set $\{\frac{k\pi}{N}:0\le k\le
2N-1\}$. The following is immediate from
Lemma~~\ref{noncom_cover_lem}.

\begin{cor}   \label{nonc_ratio_cov_cor}
Let $\vp:R\to S$ be a tame covering of polygonal surfaces. Then
one of them is rational iff the other is as well.
\end{cor}

A {\em (tame) covering of translation surfaces} is a (tame) covering of
polygonal surfaces $\vp:R\to S$, where $R$ and $S$ are translation
surfaces.

\begin{cor}   \label{noncom_trans_cov_cor}
A rational polygonal surface $P$ has a unique minimal tame
covering $\vp_t:S(P)\to P$ by a translation surface. The minimality
means that if $\vp:R'\to P$ is any tame covering by a translation
surface then there is a tame covering of translation surfaces
$\psi:R'\to S(P)$ such that $\vp=\vp_t\circ\psi$.
\begin{proof}
Follow the proof of Corollary~~\ref{trans_cov_cor} and use
Proposition~~\ref{noncom_cover_prop}.
\end{proof}
\end{cor}

The surface $S(P)$ in Corollary~~\ref{noncom_trans_cov_cor} $S(P)$ is the
(canonical) translation surface of $P$ and $\vp_t:S(P)\to P$ is
the {\em canonical translation covering}. Its degree is
equal to $|\hol_r(P)|$.

%\medskip

%We will illustrate the preceding material with examples.
%
%
\begin{exa}      \label{holo_cover_exa}
{\em Let $P$ be a noncompact, rational polygonal surface with
boundary, and let $DP$ be its doubling. Then
$\hol_r(DP)=\hol_r(P)\cap\iso_0(\RR)$ has index $2$ in
$\hol_r(P)$. We have $S(DP)=S(P)$. The canonical translation
covering $\vp_t:S(P)\to P$ is the composition of the canonical
translation covering $\vp_t:S(P)\to DP$ and the $2$-to-$1$
projection $\psi:DP\to P$.

}
\end{exa}

\section{Geodesic flow and directional flows}         \label{flows}
Let $P$ be a polygonal surface, let $\Si\subset P$ be the singular
set, and let $\bo P$ be the boundary. The {\em tangent plane}
$T_{\xi}P$ is defined for $\xi\in\inter(P)\setminus\Si$. For
$\xi\in\bo P\setminus\Si$ we denote by $T_{\xi}P$ the {\em tangent
half-plane}, i. e., the quotient of $\RR$ by the orthogonal
reflection about the side of $P$ containing $\xi$. For
$\xi\in\Si\cap\inter(P)$ (resp. $\xi\in\Si\cap\bo P$) we denote by
$T_{\xi}P$ the {\em tangent cone} (resp. {\em tangent wedge}). The
space $TP=\cup_{\xi\in P}T_xP$ is the tangent bundle. The
euclidean norm on $T_{\xi}P$ defines the set $U_{\xi}P\subset
T_{\xi}P$ of unit vectors; $UP=\cup_{\xi\in P}U_{\xi}P$ is the
{\em unit tangent bundle}.

Geodesic curves in $P$ are straight in local coordinates. By {\em
geodesics} we will mean these curves, parameterized by arclength.
We will use the notation $\{\xi(t):t\in\R\}$; it means
that the point mass is freely moving on $P$ with the unit speed.
Thus, $\dot{\xi}(t)\in U_{\xi(t)}P$ is the velocity vector. We
continue geodesics through points in $\bo P\setminus\Si$ using
the usual reflection law. However, continuation through $\Si$
is not defined.\footnote{It can be defined through wedge
points with wedge angles $\pi/n$ and cone points with cone angles
$2\pi/n$. Since such points rarely occur, we do not treat them
here differently from others.} Because of this, not all geodesics
are parameterized by $\R$; those that are parameterized by
half-lines and finite intervals are {\em singular geodesics}. A
geodesic curve of finite length is a {\em geodesic segment}.
Singular geodesic segments are the {\em saddle connections}; their
endpoints belong to $\Si$.

In view of this, the {\em geodesic flow} $G^t:UP\to UP$ is not
defined for all $t\in\R$ on the union
of all singular geodesics. Let $\xi\in P\setminus\Si$. The set of
$v\in U_{\xi}P$ that yield singular geodesics is countable, thus
the {\em singular set} $(UP)_{\sing}\subset UP$ has positive
codimension. Although the flow $G^t$ is defined for all $t$ on the
{\em regular set} $UP\setminus(UP)_{\sing}$ only, we will use the
notation $G^t:UP\to UP$. Thus, we will ignore
that $G^t(v)$ is defined for all $t$ only if $v\in
UP\setminus(UP)_{\sing}$. The reader should keep this in mind, and
make obvious adjustments, whenever necessary.

The reason we can ignore the singular set is that its volume with
respect to the invariant liouville measure on $UP$ is zero.
Indeed, locally $UP\sim U\times\RR$; the liouville measure is the
product of respective lebesgue measures. Thus, sets of smaller
dimension have measure zero. We will denote by $d\la_P, d\la_U$
the lebesgue measures on $P,U$ respectively. For the liouville
measure we have $d\mu=d\la_Pd\la_U$. Note that $d\la_P, d\la_U$
and hence $d\mu$ are determined up to a factor. If
$\vol(P)<\infty$, then we can normalize so that all relevant measures
have volume one.

\subsection{Rational polygonal surfaces and directional geodesic
flows}   \label{direction_sub}
%
%\hfill \break
If $O\subset P$ is a coordinate neighborhood, we
denote by $UO\subset UP$ the set of vectors with base points in
$O$. The $(\iso(\RR),\RR)$-structure of $P$ induces
representations $UO=O\times U$, and hence projections $p_O:UO\to
U$. If $\het=p_O(\xi)$, then $\het$ is the ``direction'' of the
vector $\xi\in UP$. The maps $p_O$ are defined up to the action of
$O(2)$ on $U$. If $g\in\iso(\RR)$, let $\bg\in O(2)$ be its linear
part. Choosing a particular $p_O$ for some $O\subset P$, continuing
the projection along a closed path, say $\ga$, and returning to
$O$, we obtain the projection $p'_O=\overline{\ho(\ga)}\circ p_O$.
Thus, the direction of a tangent vector is defined modulo the
action of $\hol_r(P)$ on the unit circle.

If $P$ be rational, then the quotient space $U/\hol_r(P)$ exists. Thus, the local projections
$p_O$ define the mapping $p:UP\to U/\hol_r(P)$. Let $\bet\in
U/\hol_r(P)$. Since the flow $G^t$ commutes with the local
projections, the set $UP_{\bet}=\{\xi\in UP:\bet(\xi)=\bet\}$ is
$G^t$-invariant. Varying $\bet\in U/\hol_r(P)$, we obtain a
measurable decomposition $UP=\cup_{\bet\in U/\hol_r(P)}UP_{\bet}$
into $G^t$-invariant sets. We will use the following convention.
Let $\het\in U$ and let $\bet\in U/\hol_r(P)$ be its projection.
We set $UP_{\het}=UP_{\bet}$. We denote by $G^t_{\het}$ the
restriction of the geodesic flow to $UP_{\het}$. The liouville
measure on $UP$ decomposes into {\em directional liouville
measures}: $d\mu=\int_{\het\in U/\hol_r(P)}d\mu_{\het}$. We
summarize the preceding discussion as follows.
\begin{prop}    \label{direction_lem}
Let $P$ be a rational polygonal surface. Then the geodesic flow
$(UP,G^t,d\mu)$ decomposes:
\begin{equation}   \label{direction_eq}
(UP,G^t,d\mu)=\int_{\het\in
U/\hol_r(P)}(UP_{\het},G^t_{\het},d\mu_{\het})d\het.
\end{equation}

The decomposition is compatible with tame coverings of polygonal
surfaces.
\end{prop}

We will call $(UP_{\het},G^t_{\het},d\mu_{\het})$ the {\em
directional (geodesic) flows} on the rational polygonal surface $P$.
Let $\het\in U$. Then $UP_{\het}=\{v\in
UP:\,p(v)\in\hol_r(P)\cdot\het\}$. Let $q:UP\to P$ be the
projection. Then $q:UP_{\het}\to P$ is onto; for any point $\xi\in
P$ the fibre $q^{-1}(\xi)$ consists of tangent vectors $v\in
U_{\xi}P$ whose directions belong to the finite orbit
$\hol_r(P)\cdot\het$.

We specialize this discussion to a translation surface, say $S$.
Then $US_{\het}$ consists of unit tangent vectors $v\in US$ with
direction $\het$. Identifying this set with $S$, we obtain the
{\em linear flow $L_{\het}^t$ on $S$ in direction} $\het$. The
orbits of $L_{\het}^t$ are the geodesics $\xi(t)$ on $S$ with
direction $\het$. The flows $L_{\het}^t$ preserve the lebesgue
measure $\la_S$.

\begin{cor}   \label{direction_cor}
Let $P$ be a rational polygonal surface. Let $S=S(P)$ be the
corresponding translation surface, and let $\vp_t:S\to P$ be the
minimal tame covering. Let $\het\in U$ be a direction whose isotropy subgroup in $\hol_r(P)$ is trivial.
The projection $\vp_t:S\to P$
induces the isomorphism of flows
$(S,L_{\het}^t,\la_S)\simeq(UP_{\het},G^t_{\het},d\mu_{\het})$.
\begin{proof}
We have already proved the claim in the special case of
translation surfaces. The general case follows from the canonical covering
$\vp_t:S\to P$, by Proposition~~\ref{direction_lem}.
\end{proof}
\end{cor}

\begin{rem}     \label{sing_dir_rem}
{\em Let $\het\in U$ be such that its isotropy subgroup in $\hol_r(P)$ is nontrivial.
We will say that $\het$ is a singular direction. A rational boundaryless polygonal
surface has no singular directions. A rational polygonal
surface with a  boundary has a finite nuber of them. Let $P$ be a
rational polygonal surface with a  boundary; let $\het\in U$ be a singular direction.
Then the relationship between $(S,L_{\het}^t,\la_S)$ and $(UP_{\het},G^t_{\het},d\mu_{\het})$
is not an isomorphism, as in Corollare~~\ref{direction_cor}. Instead, it is a $2$-to-$1$
covering of flows $q:(S,L_{\het}^t,\la_S)\to(UP_{\het},G^t_{\het},d\mu_{\het})$
This is contained in \cite{Gut84} for compact polygonal surfaces. The discussion in \cite{Gut84}
extends mutatis mutandis to all rational  polygonal surfaces.

}
\end{rem}

%\medskip

%
\subsection{Arithmetic polygonal surfaces}   \label{arithm_sub}
\hfill \break The simplest compact boundaryless polygonal surfaces
are flat tori. These are translation surfaces $T=\RR/L$, where
$L\subset\RR$ is a lattice, i. e., a closed cocompact subgroup.
The integer lattice $\ZZ\subset\RR$
yields the {\em standard torus} $T_0=\RR/\ZZ$. A compact
translation surface (of positive genus) is {\em arithmetic} if it
admits a translation covering $\vp:S\to T$ whose branching locus
is a point. Modifying the translation structure by $\gl$, if need
be, we can assume that our translation covering is $\vp:S\to T_0$
and the branching locus is $\{0\}/\ZZ$. Representing $T_0$ by
the unit square, we obtain a unit square tiling of $S$. This
explains the term {\em square tiled} for arithmetic translation
surfaces. Another popular name for these surfaces is {\em
origami} \cite{Sm04}. There are several characterizations of (compact)
arithmetic translation surfaces \cite{Gut84,GJ96,GJ00}. A {\em
compact, rational polygonal surface} $P$ is {\em arithmetic} if
its translation surface $S(P)$ is arithmetic.\footnote{This definition
is already in the literature, if only implicitly.}  We will extend these
notions to noncompact polygonal surfaces.

\begin{defin}      \label{noncom_arith_def}
{\rm Let $P$ be a rational polygonal surface. Then $P$ is {\em
arithmetic} if its translation surface $S(P)$ admits a tame
translation covering of a flat torus whose branching locus is (at
most) a point.

}
\end{defin}

Since we view polygons as polygonal surfaces,
Definition~~\ref{noncom_arith_def} applies to them as well. We
will give a few examples of noncompact arithmetic polygons.

\begin{exa}      \label{arith_polyg_exa}
{\em 1. We will use the notation of
Example~~\ref{period_obst_exa}. Let
$\Om=B_0\setminus\{\cup_{k\in\Z}(R(a,b;\xi,\eta)+(k,0))\}$ be the standard
infinite band with a periodic sequence of rectangular obstacles.
See figure~~\ref{rect_loren_band}. Then $\Om$ is arithmetic iff
$\eta,a,b\in\QQ$. 2. We will use the notation of
Example~~\ref{doubl_period_obst_exa}. Let
$\Om=\RR\setminus\{\cup_{(m,n)\in\ZZ}(R(a,b;\xi,\eta)+(m,n))\}$ be
the plane with a $\ZZ$-periodic configuration of rectangular
obstacles. Then $\Om$ is arithmetic
iff $a,b\in\QQ$.

}
\end{exa}

\begin{figure}[htbp]
\begin{center}
\input{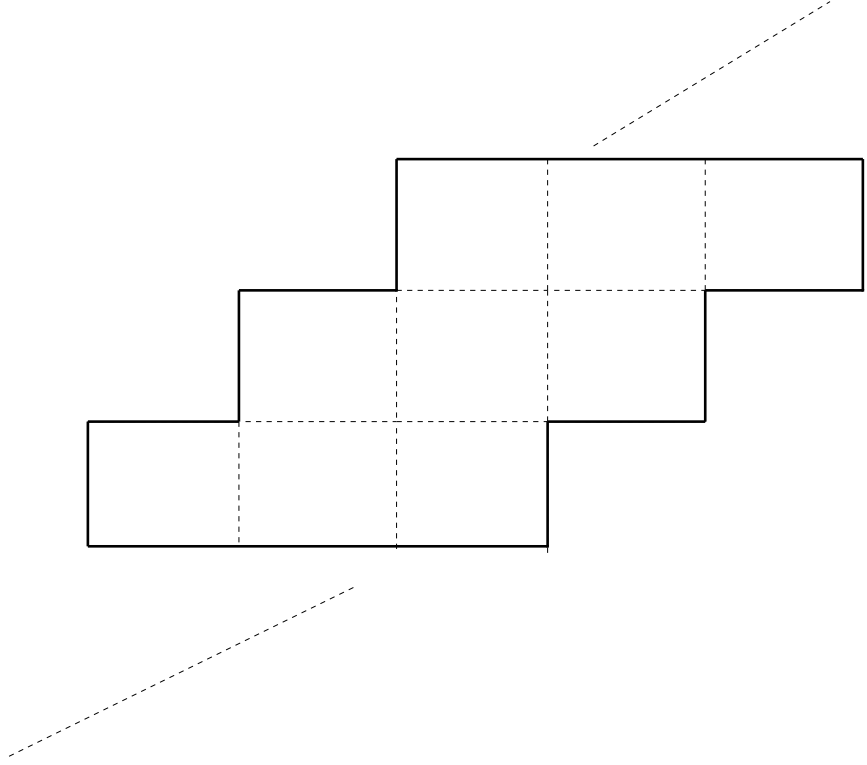}
\caption{A noncompact origami translation surface.}
\label{squar_tiled}
\end{center}
\end{figure}

\begin{exa}      \label{origami_exa}
{\em Let $P\subset\RR$ be a polygon. We say that $P$ is {\em
drawn on the integer lattice} if there is a set $I\subset\ZZ$ such that
$P=\cup_{(m,n)\in I\subset\ZZ}R_{(m,n)}$. In particular, $P$ is an arithmetic polygon;
it is compact iff $|I|<\infty$. Each side of $P$ is a horizontal or a vertical
segment with integer endpoints. Gluing some of the sides by integer translations,
we obtain an arithmetic polygonal surface. If no sides of $P$ remain unglued,
we obtain an arithmetic translation surface.

We may call these polygonal surfaces {\em square tiled}. Any arithmetic translation surface
is $\mac$-equivalent to a square tiled surface. A square tiled polygon $P$ is naturally partitioned
by horizontal rows (resp. vertical columns). A row (resp. column) is the union of a maximal
set of horizontally (resp. vertically) adjacent squares. This suggests a particular
translation surface $S$ represented by $P$. To obtain $S$, we glue the left (resp. lower) side
of every row (resp. column) with its right (resp. upper) side. Following \cite{Sm04}, we call $S$
the {\em origami translation surface} corresponding to $P$.

Figure~~\ref{squar_tiled} shows a particular noncompact origami surface $S$.
Note that each row and each column of $S$ consists of $3$ unit squares.
It is invariant under the translation by $(1,1)$ and
has the shape of an infinite stairway going from South-East to North-West.
Recurrence of geodesics on such stairways is studied in \cite{HW09}.
See also section~~\ref{transi}.

}
\end{exa}

Directional flows $G_{\het}^t$ on a compact, arithmetic polygonal
surface $P$ satisfy a dichotomy, depending on whether $\het$ is a
{\em rational direction} or an {\em irrational direction} (for $P$). In order to define this notion,
we assume first that $S(P)$ is square tiled. We use the
cartesian coordinates to identify the circle $U$ of directions
with $\R/2\pi\Z$; the {\em slope of a direction} is $\tan\het$.
Then $\het$ is {\em rational} (resp. {\em irrational}) if
$\tan\het\in\QQ\cup\{\infty\}$ (resp. $\tan\het\notin\QQ\cup\{\infty\}$). Let
now $P$ be any compact, arithmetic polygonal surface, and let
$S=S(P)$. Let $g\in\gl$ be such that $S_1=g\cdot S$ is square
tiled. A direction $\het\in U$ is rational for $P$ if $g\cdot\het$
has a rational slope. The set of slopes of rational directions for
$P$ does not depend on the choices involved. It has the
form $g^{-1}\cdot(\QQ\cup\{\infty\})$, where $g^{-1}\in\gl$ is a fractional
linear transformation, hence countable.

%A geodesic segment on $P$ whose endpoints belong to the singular
%set is a {\em saddle connection}.

%
\begin{thm}    \label{gut_dichot_thm}
{\em (See \cite{Gut84}.)} Let $P$ be a compact, arithmetic
polygonal surface. If $\het$ is irrational, then $G_{\het}^t$ is
uniquely ergodic. If $\het$ is rational, then every geodesic in
direction $\het$ is either periodic or a saddle connection.
\end{thm}

Note that the above definition of rational and irrational
directions applies to noncompact, arithmetic polygonal surfaces
as well. We conjecture that the  dichotomy of
Theorem~~\ref{gut_dichot_thm} extends to arbitrary arithmetic
polygonal surfaces. We point out that the noncompact
version of Theorem~~\ref{gut_dichot_thm} should take into account
the possibility of {\em transient geodesics}. We discuss this in
the next section.

\section{Conservativeness and dissipation for noncompact
polygonal surfaces}   \label{transi}
First, we will recall  basic notions pertaining to recurrence and
transience in dynamical systems \cite{Aa97,Sc77}. For simplicity
of exposition, we gear the discussion to the dynamics with time
$\Z$. We will use the shorthand {\em $\nu$-a.e.} to mean {\em almost every with respect to (the measure) $\nu$}.
The dynamical system $(X,\tau,\nu)$ is {\em conservative} if
for every measurable set $B\subset X$  and for $\nu$-a.e. point
$x\in B$ there is $n=n(x)>0$ such that $\tau^n x \in B$. It is
{\em dissipative} if there is a measurable subset $A\subset X$
such that $X=\cup_{k\in\Z}\tau^{k}A$ and the sets
$\tau^{k}A\cap\tau^{l}A=\emptyset$ for $k\ne l$. Every dynamical
system uniquely decomposes as a disjoint union of its {\em
conservative part} and the {\em dissipative part}: $X=X_{\con}\cup
X_{\dis}$.  Points $x\in X_{\con}$ (resp. $x\in X_{\dis}$) are the
{\em recurrent} (resp. {\em transient}) points in $X$. If
$\nu(X)<\infty$, then, by the poincare recurrence theorem,
$\nu$-a.e. point is recurrent; equivalently, $X=X_{\con}$. If
$\nu(X)=\infty$, which is our case, the dissipative part may be
nontrivial. Note that $\nu(X_{\dis})>0$ iff
$\nu(X_{\dis})=\infty$.

It is straightforward to reformulate the above notions for the
dynamics with time $\R$, i. e., for measure preserving flows. We
will study the recurrence of geodesic flows on certain noncompact
polygonal surfaces. The following observation will allow us to
replace the dynamics with time $\R$ by the dynamics with time
$\Z$.

\begin{prop}      \label{cros_sect_prop}
Let $(Y,T^t,\mu)$ be a flow, let $X\subset Y$ be a cross-section,
and let $(X,\tau,\nu)$ be the induced return transformation. Then
the time $\R$ dynamical system $(Y,T^t,\mu)$ is conservative iff
the time $\Z$ dynamical system $(X,\tau,\nu)$ is conservative.
\end{prop}

\subsection{Geodesic flows on $\Z$-periodic polygonal surfaces}    \label{periodic_sub}
\hfill \break We begin by defining the notion.

\begin{defin}    \label{per_surf_def}
Let $\tP$ be a noncompact polygonal surface, and let $\Ga$ be a
countably infinite group acting by
isometries on $\tP$. Suppose that $\Ga$ acts freely and cocompactly.
Then $\tP$ is a {\em $\Ga$-periodic polygonal surface}.
\end{defin}

In our examples, $\Ga=\Z$ or $\Ga=\ZZ$. For instance, the surfaces
in Example~~\ref{period_obst_exa} and Figure~~\ref{inf_cyl_wind}
are $\Z$-periodic; the surface in
Example~~\ref{doubl_period_obst_exa} is $\ZZ$-periodic. If $\Ga$
is not specified, we will speak of {\em periodic polygonal
surfaces}.

Let $\tP$ be a $\Ga$-periodic polygonal surface. Then $P=\tP/\Ga$
is a compact polygonal surface; let $p:\tP\to P$ be the
projection. Let $U\tP,UP$ be the unit tangent bundles for $\tP,P$;
let $\tG^t,G^t$ be the respective geodesic flows; let $\tmu,\mu$
be the liouville measures for $U\tP,UP$ respectively. The action
of $\Ga$ on $\tP$ uniquely extends to a free, cocompact action on
$U\tP$. We have $UP=U\tP/\Ga$; let $q:U\tP\to UP$ be the
projection. The polygonal surfaces $\tP$ and $P$ are rational or
not rational simultaneously. Suppose they are rational; let
$(U\tP_{\het},\tG_{\het}^t,\tmu_{\het})$ and
$(UP_{\het},G_{\het}^t,\mu_{\het})$ be the respective directional
geodesic flows. By Proposition~~\ref{direction_lem}, the
projection $q:U\tP\to UP$ is compatible with the decomposition
equation~~\eqref{direction_eq}, yielding the directional
projections
$q_{\het}:(U\tP_{\het},\tG^t_{\het},\tmu_{\het})\to(UP_{\het},G^t_{\het},\mu_{\het})$.
Note that $q$ and all $q_{\het}$ are coverings of flows.

The geodesic flow $(U\tP,\tG^t,\tmu)$ (resp. directional flow
$(U\tP_{\het},\tG^t_{\het},\tmu_{\het})$) is a  {\em skew product
over the flow $(UP,G^t,\mu)$ (resp.
$(UP_{\het},G^t_{\het},\mu_{\het})$) with the fibre $\Ga$.} To
simplify the exposition, we recall the concept of skew products
for time $\Z$ dynamical systems. The  time $\R$ case  is similar.
See \cite{Co09} for more information.

Let $(X,\tau,\nu)$ be a measure preserving automorphism. Let $\Ga$
be a countably infinite abelian group, and let $\nu_{\Ga}$ be the counting
measure.\footnote{The notion of skew products makes sense for any
locally compact group. We restrict the exposition to this class of groups because
in our applications, $\Ga=\Z$ or $\Ga=\ZZ$.}
Let $\vp:X\to\Ga$ be a measurable mapping. Set
$\tX=X\times\Ga$, $\tnu=\nu\times\nu_{\Ga}$, and
\begin{equation}     \label{produit_gauche_eq}
\ttau(x,g) = (\tau x, g+\vp(x)).
\end{equation}
The dynamical system $(\tX,\ttau,\tnu)$ is the {\em skew product over
$(X,\tau,\nu)$ with the fibre $\Ga$ and the displacement function $\vp$}.

Let $(\tY,\tT^t,\tmu)$ be a measure preserving flow which is a
skew product with fibre $\Ga$ over the flow $(Y,T^t,\mu)$. Let
$q:(\tY,\tT^t,\tmu)\to(Y,T^t,\mu)$ be the projection. Let
$X\subset Y$ be a cross-section for $(Y,T^t,\mu)$. Let $\nu$ be
the induced measure on $X$, and let $(X,\tau,\nu)$ be the poincare
map. Set $\tX=q^{-1}(X)\subset\tY$. Then $\tX$ is a cross-section
for $(\tY,\tT^t,\tmu)$. Let $(\tX,\ttau,\tnu)$ be the
corresponding poincare map. Then
$(\tY,\tmu)=(Y\times\Ga,\mu\times\nu_{\Ga})$ as measure spaces.
This induces the isomorphism
$(\tX,\tnu)=(X\times\Ga,\nu\times\nu_{\Ga})$ and determines a
mapping $\vp:X\to\Ga$ such that $\ttau:\tX\to\tX$ is given by
equation~~\eqref{produit_gauche_eq}. Thus, $(\tX,\ttau,\tnu)$ is
the skew product over $(X,\tau,\nu)$; the displacement function
$\vp:X\to\Ga$ is uniquely determined by
$q:(\tY,\tT^t,\tmu)\to(Y,T^t,\mu)$ and the choice of cross-section
$X\subset Y$. The following example illustrates this material in
the context of periodic polygonal surfaces.
\begin{exa}    \label{cylin_exa}
{\em We will use two subgroups $\Z\subset\RR$: The {\em horizontal
$\Z_h=\{(k,0)$ and the vertical $\Z_v=\{(0,k)\}$.} The quotient
$\tS=\RR/\Z_h$ is a flat vertical cylinder; $\Z_v$ acts on $\tS$
by vertical translations. The action is free and cocompact. The
quotient $S=\tS/\Z_v$ is the flat unit torus. We view $S$ as the
unit square $R=\{(x,y):0\le x,y \le 1\}$ with its opposite sides
glued in the usual way. Let $U\tS,US$ be the unit tangent bundles
of $\tS,S$; let $(U\tS,\tG^t,\tmu)$ and $(US,G^t,\mu)$ be the
respective geodesic flows. Identifying the {\em circle of
directions} with $[0,2\pi)$, we have $US=\{(x,y,\het):(x,y)\in
R,\het\in[0,2\pi)\}$. The set $X=\{(x,0,\het)\}$ is a
cross-section for $(US,G^t,\mu)$. With the notation $(x,\het)\in
X$, the poincare map is given by
$\tau(x,\het)=((x+\cot\het)\mod1,\het)$. Let $p:\tS\to S$ and
$q:U\tS\to US$ be the projections. We identify $U\tS$ with
$\{(x,y,\het):(x,\het)\in X,y\in\R\}$. Then
$\tX=q^{-1}(X)=\{(x,k,\het):(x,\het)\in X,k\in\Z\}$. Thus,
$\tX=X\times\Z$ up to a set of measure zero. Set $\vp(x,\het)=1$
if $0<\het<\pi$, $\vp(x,\het)=-1$ if $\pi<\het<2\pi$. The poincare
map for the cross-section $\tX$ is given by
$\ttau((x,\het),k)=(\tau(x,\het),k+\vp(x,\het))$. Thus,
$\vp:X\to\Z$ is the displacement function.

In this example $\tS,S$ are translation surfaces; $\het$
corresponds to the direction of a tangent vector. Fixing $\het$,
we obtain the subset $X_{\het}\subset X$ which is a cross-section
for the directional flow $(US_{\het},G^t_{\het},\mu_{\het})$. We
have $X_{\het}=\R/\Z$. The directional poincare map is given by
$\tau_{\het}\cdot x=(x+\cot\het)\mod1$. Set
$q^{-1}(X_{\het})=\tX_{\het}=(\R/\Z)\times\Z$. The directional
poincare map for $\tS$ is
$\ttau_{\het}\cdot(x,k)=((x+\cot\het)\mod1,k+\vp_{\het}(x))$. The
{\em directional displacement functions} $\vp_{\het}:\R/\Z\to\Z$
are $\vp_{\het}=1$ for $0<\het<\pi$, and $\vp_{\het}=-1$ for
$\pi<\het<2\pi$.

}
\end{exa}

Example~~\ref{cylin_exa} shows  the importance of choice of
cross-sections for (directional) geodesic flows. For polygonal
surfaces with boundaries, there is a canonical choice: The set of
tangent vectors with base points on the boundary. We will now give
the relevant definitions.

Let $P$ be an arbitrary polygonal surface. Let $A\subset P$ be a
closed subset. We denote by $U_AP\subset UP$ the set of tangent
vectors $v\in UP$ whose base points belong to $A$. Suppose that
$P$ has a boundary, $\bo P$. Let $BUP\subset UP$ be the smallest
$G^t$-invariant set containing $U_{\bo P}P$. Let
$\mu_B=\mu|_{BUP}$. By definition, $U_{\bo P}P$ is a cross-section
for the flow $(BUP,G^t,\mu_B)$. Let $\nu$ be the induced measure
on $U_{\bo P}P$. We will use the following terminology: The flow
$(BUP,G^t,\mu_B)$ is the {\em billiard flow} for the polygonal
surface $P$; the set $U_{\bo P}P\subset BUP$ is the {\em standard
cross-section} for the billiard flow; the induced transformation
$(U_{\bo P}P,\tau,\nu)$ is the {\em billiard map} for $P$.

For $N\in\N$ denote by $D_N\subset O(2)$ the group generated by
orthogonal reflections about two axes forming the angle $\pi/N$.
It is the {\em dihedral group} of order $2N$. Suppose now that $P$
is a rational polygonal surface. Then $\hol_r(P)=D_N$ where
$N=N(P)$. If $P$ is a rational polygon, then $N$ is the least
common denominator of its angles. We identify $U/D_N$ and
$[0,\pi/N]$. For $\het\in[0,\pi/N]$ set $U_{\bo
P}P_{\het}=UP_{\het}\cap U_{\bo P}P$, $BUP_{\het}=BUP\cap
UP_{\het}$.

\begin{defin}      \label{transvers_def}
Let $P$ be a rational polygonal surface with a boundary. A
direction $\het\in[0,\pi/N]$ is {\em transversal to the boundary
of $P$} if $\mu_{\het}(BUP_{\het})>0$.
\end{defin}

If $\het$ is transversal to $\bo P$, we denote by $b\mu_{\het}$
the restriction of $\mu_{\het}$ to $BUP_{\het}$; then
$(BUP_{\het},G^t_{\het},b\mu_{\het})$ is the {\em billiard flow in
direction $\het$}. The set $U_{\bo P}P_{\het}\subset BUP_{\het}$
is the {\em standard directional cross-section}. Let $\nu_{\het}$
be the induced measure on $U_{\bo P}P_{\het}$. The induced
dynamical system $(U_{\bo P}P_{\het},\tau_{\het},\nu_{\het})$ is
the {\em directional billiard map}. When $P\subset\RR$ is a
compact rational polygon, this is the standard terminology
\cite{Gut03}. When the surface $P$ is clear from the context, we
will simply speak of {\em transversal directions}. We will now
illustrate the above material with examples.

Let $S$ be a translation surface. Let $O\subset S$ be a polygon,
in general disconnected, such that $S\setminus O$ is connected.
The polygonal surface $P=S\setminus\inter(O)$ is a {\em
translation surface with polygonal obstacles}. A {\em barrier}
$L\subset S$ is a connected polygon without interior, e. g., a
linear segment. Making the billiard ball bounce off of each side
of $L$, we obtain the polygonal surface which is a {\em
translation surface with a barrier}. We will denote it by
$P=S\setminus L$. Combining the two notions, we come to the
concept of {\em translation surfaces with polygonal obstacles
and/or barriers}. Note that these are polygonal surfaces with
boundaries.

\begin{figure}[htbp]
\begin{center}
\input{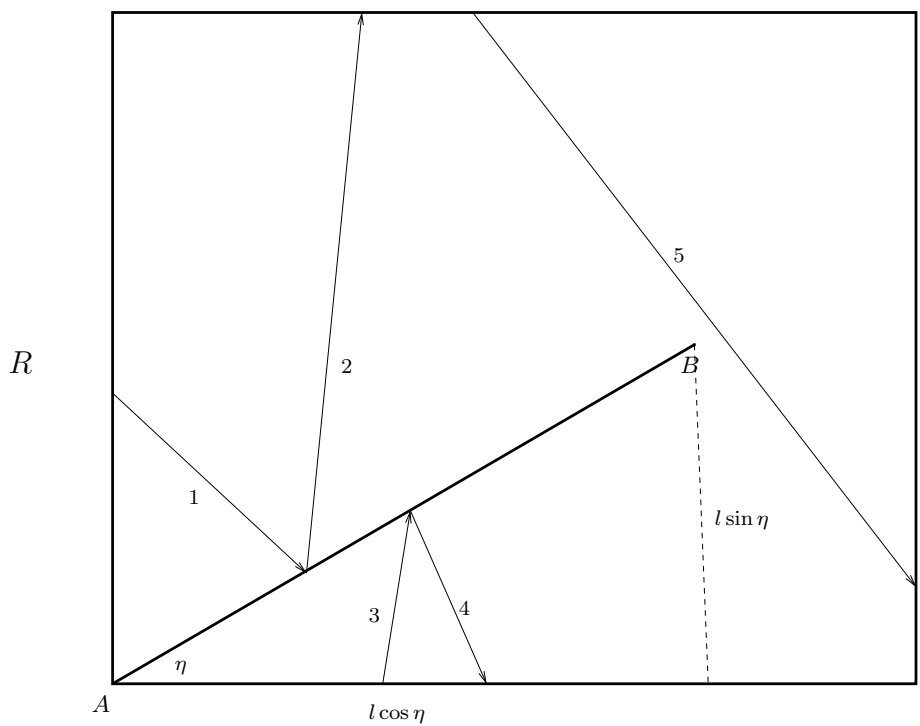}
\caption{A billiard orbit in the standard torus with a linear
barrier.} \label{tor_with_bar}
\end{center}
\end{figure}

\begin{exa}     \label{tor_with_bar_exa}
{\em Let $S=\RR/\ZZ$ be the flat unit torus. A {\em linear
barrier} $L\subset S$ is determined by one of its end points, say
$A\in S$, its direction $\eta$, and its length $l$. We assume
without loss of generality, that $A=(0,0)$ and that $0\le\eta
<\pi/2$. The other endpoint of the barrier is
$B=(l\cos\eta\mod1,l\sin\eta\mod1)$. We require that $B\ne A$;
otherwise $P=S\setminus L$ is disconnected. If $\eta$ is
$\pi$-rational, this yields un upper bound on $l$; otherwise $l$
is arbitrary.

As in Example~~\ref{cylin_exa}, we represent the flat unit torus
by the unit square $R$. Let $l<\frac{1}{\cos\eta}$. Then
$P=S\setminus L$ is represented by $R$ with the linear segment
$[A=(0,0),B=(l\cos\eta,l\sin\eta)]$. Figure~~\ref{tor_with_bar}
shows a billiard orbit in $P$. Denote by $P(l,\eta)$ the above
surface. We will now illustrate Definition~~\ref{transvers_def}.
The group $\hol_r(P)$ is generated by a single reflection, thus
$N=1$. Let $0\le\het\le\pi$. The set $BUP_{\het}$ is formed by the
billiard orbits crossing the barrier in direction $\het$. Up to a
universal normalizing factor,
$\mu_{\het}(BUP_{\het})=l|\sin(\het-\eta)|$ \cite{Gut86}. Thus,
the only direction which is not transversal to the boundary of
$P(l,\eta)$ is $\eta$.

}
\end{exa}

\begin{lem}      \label{transvers_lem}
Let $P$ be a compact, rational polygonal surface with a boundary.
Then all but a finite number of directions are transversal.
\begin{proof}
Let $\eta$ be the direction of a side in $\bo P$. The argument in
Example~~\ref{tor_with_bar_exa} shows that every direction
$\het\neq\eta,\eta+\pi$ is transversal.
\end{proof}
\end{lem}
\begin{rem}     \label{transvers_rem}
{\em In fact, a stronger statement holds. But for a very special
class of surfaces, every direction is transversal to $\bo P$. For
$P$ in that class, all but two directions are transversal. We will
not use these facts in what follows.

}
\end{rem}

\begin{thm}    \label{per_rat_recurr_thm}
Let $\tP$ be a $\Z$-periodic polygonal surface such that the
quotient $P=\tP/\Z$ is a translation surface with obstacles and/or
barriers. Suppose that $P$ is a rational polygonal surface and
that $N=N(P)$ is even. Then for almost every $\het\in[0,\pi/N]$
the directional flow $(U\tP_{\het},\tG^t_{\het}, {\tmu}_{\het})$
is recurrent.
\begin{proof}
Let $\het$ be an ergodic direction for $P$ which is transversal to
$\bo P$. Then $(UP_{\het},G^t_{\het},
\mu_{\het})=(BUP_{\het},G^t_{\het}, {b\mu}_{\het})$. In other
words, $\bo P$ provides a cross-section, say $X=X(\het)$, for the
directional flow $(UP_{\het},G^t_{\het},\mu_{\het})$. Let
$(X,\tau_{\het},\nu_{\het})$ be the poincare map. Let $\tX$ be the
corresponding cross-section for the directional flow
$(U\tP_{\het},\tG^t_{\het},\tmu_{\het})$ on the noncompact
surface. Let $\vp_{\het}$ be the directional displacement
function. Since $N$ is even, we have satisfies
\begin{equation}   \label{centered_eq}
\int_X\vp_{\het}\,d\nu_{\het}=0.
\end{equation}
See Lemma~~5 in
\cite{CoGu10}.\footnote{We will say that a displacement function
satisfying equation~~\eqref{centered_eq} is {\em centered}.} By
\cite{At76}, the poincare map $(\tX,\ttau_{\het},\tnu)$ is
recurrent. By Proposition~~\ref{cros_sect_prop}, the flow
$(U\tP_{\het},\tG^t_{\het}, {\tmu}_{\het})$ is recurrent as well.

By \cite{KMS}, the set of ergodic directions for a compact,
rational polygonal surface has full measure. In view of
Lemma~~\ref{transvers_lem}, almost every direction is ergodic and
transversal.
\end{proof}
\end{thm}

\begin{cor}    \label{per_rat_recurr_cor}
Let $\tP$ be a $\Z$-periodic polygonal surface satisfying the
assumptions of Theorem~~\ref{per_rat_recurr_thm}. Then the
geodesic flow on $\tP$ is conservative.
\begin{proof}
Immediate from Theorem~~\ref{per_rat_recurr_thm} and
equation~~\eqref{direction_eq}.
\end{proof}
\end{cor}

Let $\tP$ be a $\Z$-periodic polygonal surface. Let $P=\tP/\Z$ be
the quotient. Suppose that i) $P$ is a rational polygonal surface;
ii) $P$ is a translation surface with a boundary and/or obstacles;
and iii) $N(P)$ is even. Then, by
Corollary~~\ref{per_rat_recurr_cor}, the geodesic flow on $\tP$ is
conservative. In view of Example~~\ref{cylin_exa}, condition i)
alone does not ensure the conservativeness of the geodesic flow on
$\tP$. The quotient surface in Example~~\ref{cylin_exa} is
rational but has no boundary. What happens if the quotient surface
has a boundary but not rational? We will partially answer this
question.

Compact euclidean polygons form a topological space. Fixing the
type of a polygon, we obtain closed subspaces in this space. Thus,
we can ask questions about the billiard in topologically typical
polygons (of a particular type). See \cite{KZ,KMS,GK}.

We fix a compact translation surface $S_0$. Let $O\subset S_0$ be
a collection of polygons/barriers of a particular type such that
$P=S_0\setminus O$ is connected. For instance, $O$ can be an
arbitrary triangle, or a quadrilateral, or a right triangle, or a
linear segment, or a disjoint union of two triangles, etc. We
denote the type of $O$ by $\ka$. Let $\ppp=\ppp(S_0,\ka)$ be the
space of polygonal surfaces $P=S_0\setminus O$, endowed with the
natural topology. The space $\ppp$ is homeomorphic to a bounded
domain in some euclidean space. Let $\tppp=\tppp(S_0,\ka)$ be the
space of $\Z$-periodic polygonal surfaces $\tP$ such that
$\tP/\Z\in\ppp(S_0,\ka)$. The topology of $\ppp(S_0,\ka)$ induces
a topology on $\tppp(S_0,\ka)$.

Let $\ppp_1\subset\ppp$ be the subset of rational polygonal
surfaces $P=S_0\setminus O$ such that the number $N=N(P)$ is even.

\begin{defin}   \label{type_def}
We say that $\ka$ is an {\em amenable type} if $\ppp_1$ is dense
in $\ppp$.
\end{defin}

Saying that a noncompact polygonal surface $\tP$ is {\em
conservative}, we will mean that the geodesic flow
$(U\tP,\tG^t,\tmu)$ is recurrent.

\begin{thm}    \label{dens_gdel_recur_thm}
Let $\ka$ be a type of polygons/barriers in $S_0$. Let
$\ppp=\ppp(S_0,\ka)$ be the topological space of polygonal
surfaces $P=S_0\setminus O,\,O\in\ka$. Let $\tppp=\tppp(S_0,\ka)$
be the corresponding space of $\Z$-periodic polygonal surfaces.

If $\ka$ is an amenable type then the set of conservative surfaces
in $\tppp$ is a dense $G_{\de}$.
\begin{proof}
Let $\ka$ be arbitrary. It is standard to check that the set
$\tppp_{\con}(S_0,\ka)$ of conservative surfaces is a countable
intersection of open subsets in  $\tppp(S_0,\ka)$, i. e., it is a
$G_{\de}$ set. Since $\ka$ is amenable, $\tppp_{\con}(S_0,\ka)$
contains a dense subset of surfaces $\tP$ such that $P=\tP/\Z$ is
rational and $N(P)$ is even. The claim now follows from
Corollary~~\ref{per_rat_recurr_cor}.
\end{proof}
\end{thm}

Let $P_0$ be a particular compact polygonal surface. As always, we
assume that $P_0$ is connected. In our applications,
$P_0$ will be an ``elementary'' polygonal
surface.\footnote{Typical examples: i) $P_0$ is the standard torus;
ii) $P_0$ is the standard cylinder.} Choose a type $\ka$ of
obstacles/barriers $O\subset P_0$. Analogously to the preceding
discussion, we define the space $\ppp=\ppp(P_0,\ka)$ of polygonal
surfaces $P=P_0\setminus O,\,O\in\ka$, and endow it with the
natural topology. Let $\tppp=\tppp(P_0,\ka)$ be the space of
$\Z$-periodic  polygonal surfaces $\tP$ such that
$\tP/\Z\in\ppp(P_0,\ka)$. The topology of $\ppp(P_0,\ka)$ induces
a topology on $\tppp(P_0,\ka)$.

By analogy with Definition~~\ref{type_def}, we say that the pair
$(P_0,\ka)$ is {\em amenable} if i) $P_0$ is a rational polygonal
surface; ii) $\ppp(P_0,\ka)$ contains a dense subspace of rational
polygonal surfaces $P=P_0\setminus O$ such that $N(P_0\setminus
O)$ is even.

\begin{thm}    \label{mor_dens_gdel_thm}
Let $(P_0,\ka)$ be an amenable pair. Then the set of conservative
$\Z$-periodic polygonal surfaces $\tP\in\tppp(P_0,\ka)$ is a dense
$G_{\de}$.
\begin{proof}
Let $S_0=S(P_0)$ be the minimal translation covering of $P_0$; let
$\vp_0:S_0\to P_0$ be the canonical projection. See
Corollary~~\ref{trans_cov_cor}. Let $\Ga\subset O(2)$ be the
corresponding subgroup, so that $P_0=S_0/\Ga$. Then $\ka$ defines
a $\Ga$-invariant type $\la$ of obstacles/barriers $O$ in $S_0$.
The action of $\Ga$ on $S_0$ uniquely extends to the actions of
$\Ga$ on associated spaces. In particular, $\Ga$ acts on the space
$\tppp(S_0,\la)$. Pulling back by $\vp_0$ yields an isomorphism of
$\tppp(P_0,\ka)$ and $\tppp(S_0,\la)^{\Ga}$. Now the proof of
Theorem~~\ref{dens_gdel_recur_thm} applies and yields the claim.
\end{proof}
\end{thm}

We will now apply Theorem~~\ref{dens_gdel_recur_thm} and
Theorem~~\ref{mor_dens_gdel_thm} to particular families of
$\Z$-periodic polygonal surfaces.

\subsection{Examples and concluding remarks}    \label{exa_sub}
\hfill \break Let $0\le a,b\le 1$, let $0\le\het<2\pi$. Let $R$ be
the standard unit square, $R=\{(x,y):0\le x,y\le 1\}$. Let $l\ge
0$ be such that $(a+l\cos\het,b+l\sin\het)\in R$. We denote by
$O(a,b,l,\het)\subset R$ the linear segment with endpoints $(a,b)$
and $(a+l\cos\het,b+l\sin\het)$.

Let $C_0$ be the polygonal surface obtained by identifying the
vertical sides of $R$ by the parallel translation
$(\xi,\eta)\mapsto(\xi+1,\eta)$. It is a compact flat cylinder; we
will refer to $C_0$ as the {\em standard cylinder}. The surface
$C_0(a,b,l,\het)=C_0\setminus O(a,b,l,\het)$ is the {\em standard
cylinder with a linear barrier}. Let $\OO\subset\R^4$ be the
bounded domain formed by $(a,b,l,\het)$ satisfying the above
conditions. We denote by $\OO(a,b)\subset\OO$ the subset obtained
by fixing the values of the first two variables. The sets
$\OO(l),\OO(a,b,l)$, etc are similarly defined.

\begin{figure}[htbp]
\begin{center}
\input{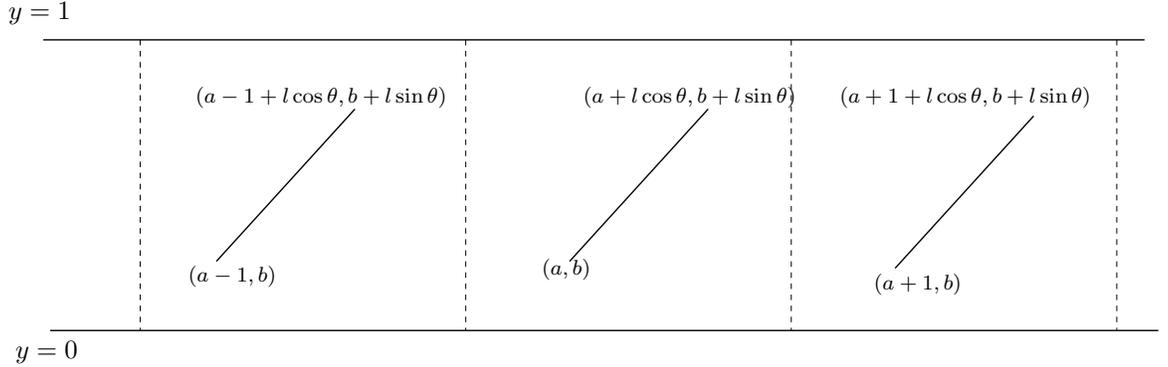}
\caption{The standard band with a $Z$-periodic collection of
tilted linear barriers.} \label{band_with_bar}
\end{center}
\end{figure}

Let $B_0$ be the standard band. See
Example~~\ref{period_obst_exa}. Set
$B_0(a,b,l,\het)=B_0\setminus\cup_{n\in\Z}\left(O(a,b,l,\het)+(n,0)\right)$ be
the standard band with a $Z$-periodic
collection of linear barriers. See figure~~\ref{band_with_bar}.
Then $B_0(a,b,l,\het)$ is a $\Z$-periodic polygonal surface, and
$B_0(a,b,l,\het)/\Z=C_0(a,b,l,\het)$.

\begin{cor}    \label{band_with_bar_cor}
1. The set of parameters in $\OO$ such that the noncompact
polygonal surface $B_0(a,b,l,\het)$ is conservative is a dense
$G_{\de}$. 2. Let $(a,b)\in R$ be arbitrary. Let
$\OO(a,b)\subset\OO$ be the corresponding subset. Then set of
parameters $(l,\het)$ such that the noncompact polygonal surface
$B_0(a,b,l,\het)$ is conservative is a dense $G_{\de}$ in
$\OO(a,b)$. 3. Let $l$ be such that the set $\OO(l)\subset\OO$ has
nonempty interior. Then set of parameters $(a,b,\het)$ such that
the noncompact polygonal surface $B_0(a,b,l,\het)$ is conservative
is a dense $G_{\de}$ in $\OO(l)$. 4. Let $(a,b,l)$ be such that
the set $\OO(a,b,l)\subset\OO$ has nonempty interior. Then set of
directions $\het$ such that the noncompact polygonal surface
$B_0(a,b,l,\het)$ is conservative is a dense $G_{\de}$ in
$\OO(a,b,l)$.
\begin{proof}
Since all claims fit into the framework of
Theorem~~\ref{mor_dens_gdel_thm}, we only need to prove that
the relevant types $(C_0,\ka)$ are amenable. In each case the
angle $\het$ is free to vary in certain intervals.\footnote{Note
that the intervals may depend on relevant parameters.} The group
in question is generated by the reflection about the horizontal
axis and about the $\het$-axis. It is finite iff
$\frac{\het}{\pi}=\frac{m}{n}$, and then $N(P)=n$. All claims now
follow from the observation that rational numbers with even
denominators are dense in any interval.
\end{proof}
\end{cor}

\begin{figure}[htbp]
\begin{center}
\input{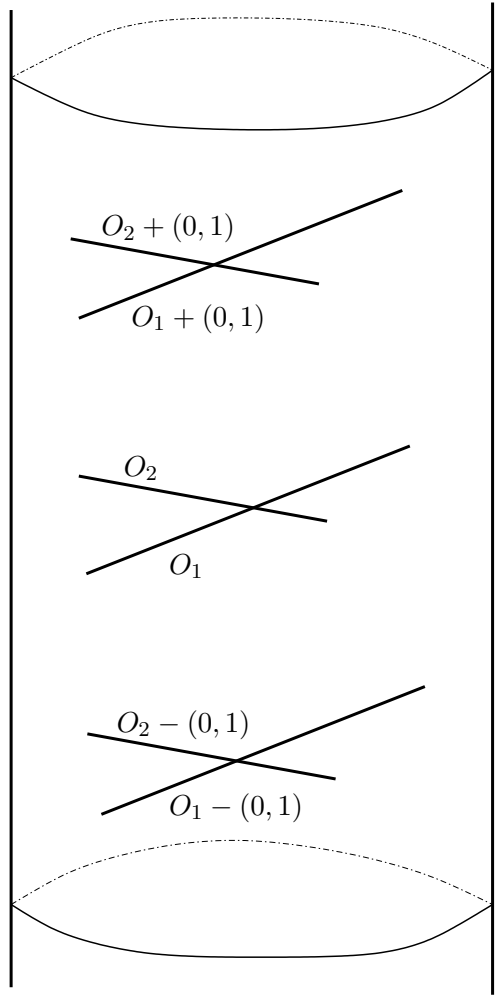}
\caption{Infinite cylinder with two $\Z$-periodic collections of
linear barriers.} \label{cyl_with_bar}
\end{center}
\end{figure}

Let $\tC_0$ be the translation surface obtained by gluing the
boundary components of $B_0$ via the parallel translation
$(x,0)\mapsto(x,1)$. The infinite flat cylinder $\tC_0$ satisfies
$\tC_0=\RR/\Z_v$. See Example~~\ref{cylin_exa}. Let $T_0=\RR/\ZZ$ be
the standard torus. Let $0\le\het_1,\het_2<2\pi$ be arbitrary
angles. Let $a_i,b_i,l_i,i=1,2$ be parameters satisfying the
conditions described prior to Corollary~~\ref{band_with_bar_cor}.
Denote by $O_i=O(a_i,b_i,l_i,\het_i)\subset R,i=1,2$ the
corresponding linear segments. The compact polygonal surface
$P=P(a_i,b_i,l_i,\het_i:i=1,2)=T_0\setminus\left(O_1\cup
O_2\right)$ is the {\em torus with two linear barriers}. If
$\het_2-\het_1\ne k\pi$, the barriers are transversal. Otherwise,
they may overlap, forming a single barrier. Denote by
$\tP=\tP(a_i,b_i,l_i,\het_i:i=1,2)=C_0\setminus\left((\Z\cdot O_1)\cup
(\Z\cdot O_2)\right)$ the $\Z$-periodic covering surface. Then $\tP$
is the {\em infinite cylinder with two $\Z$-periodic collections
of linear barriers}. See figure~~\ref{cyl_with_bar}.

Let $\OO\subset\R^4$ be the bounded domain defined prior to
Corollary~~\ref{band_with_bar_cor}. For $r_i\in\OO,i=1,2$ let
$\tP(r_1,r_2)$ be the corresponding infinite cylinder with two
$\Z$-periodic collections of linear barriers. This is a family of
$\Z$-periodic polygonal surfaces parameterized by points in
$\OO\times\OO\subset\R^8$.

\begin{cor}    \label{cyl_with_bar_cor}
The set of parameters $(r_1,r_2)\in\OO\times\OO$ such that the
noncompact polygonal surface $\tP(r_1,r_2)$ is conservative is a
dense $G_{\de}$.
\begin{proof}
Set $P(r_1,r_2)=\tP(r_1,r_2)/\Z=T_0\setminus\left(O_1\cup
O_2\right)$. The claim fits into the framework of
Theorem~~\ref{dens_gdel_recur_thm}. Hence, it suffices to prove that the
type $(T_0,\ka)$ of two linear barriers in the torus is amenable.
If we fix non-angular parameters, the angles $\het_1,\het_2$ are
free to vary in certain intervals which in general depend on these
parameters. The surface $P(r_1,r_2)$ is rational iff
$\frac{\het_1-\het_2}{\pi}=\frac{m}{n}$,\footnote{Where $m$ and $n$
are relatively prime.} and $N(P(r_1,r_2))=n$.
The claim now follows from the observation that the condition
$\frac{\het_1-\het_2}{\pi}=\frac{m}{n}$ with $n$ even is dense in
the product of  any intervals.
\end{proof}
\end{cor}

\begin{rem}   \label{cyl_with_bar_rem}
{\em Corollary~~\ref{cyl_with_bar_cor} is the analogue of claim 1
in Corollary~~\ref{band_with_bar_cor}. Fixing some of the
coordinates in $(r_1,r_2)$, we obtain analogs of the other claims
in Corollary~~\ref{band_with_bar_cor}. The dense $G_{\de}$ claim
will hold as long as at least one of the angles $\het_i$ is free
to vary in an interval. We leave the details to the reader.

}
\end{rem}

\begin{figure}[htbp]
\begin{center}
\input{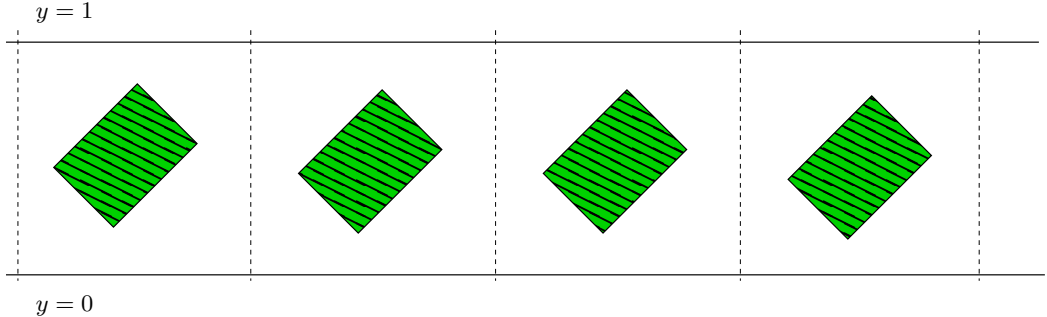}
\caption{Infinite band with a $\Z$-periodic collection of tilted
rectangular obstacles.} \label{cyl_with_rect}
\end{center}
\end{figure}

Let $R(a,b;\xi,\eta)$ be the $a\times b$ rectangle with the lower
left corner $(\xi,\eta)$. See Example~~\ref{period_obst_exa} and
Example~~\ref{qua_perio_obst_exa} for notation. Let
$\rho_{\het}R(a,b;\xi,\eta)$ be the rectangle $R(a,b;\xi,\eta)$
rotated by the angle $\het$ about its center-point. Let $R$ be the
unit rectangle. We assume that the parameters $a,b;\xi,\eta$ are
such that $\rho_{\het}R(a,b;\xi,\eta)\subset R$ for $0\le\het\le
2\pi$. We fix $a,b;\xi,\eta$ and set
$O_{\het}=\rho_{\het}R(a,b;\xi,\eta)\subset R$. Representing the
standard cylinder $C_0$ by $R$, we have $O_{\het}\subset C_0$. Set
$P_{\het}=C_0\setminus O_{\het}$. Let $\ka$ be the type of
obstacles $O_{\het}\subset C_0,0\le\het\le 2\pi$. Thus,
$\ppp=\ppp(C_0,\ka)=\{P_{\het}:0\le\het\le 2\pi\}$ is a space of
flat tori with rectangular obstacles. Note that topologically the space $\ppp$
is the circle. Let $\tppp=\tppp(C_0,\ka)$ be the
corresponding space of $\Z$-periodic polygonal surfaces. Let $B_0$
be the standard infinite band. The surface $\tP_{\het}\in\tppp$ is
$B_0$ with a $\Z$-periodic collection of tilted rectangular
obstacles $O_{\het}$. We have
$\tP_{\het}=B_0\setminus\left(\cup_{k\in\Z}(O_{\het}+(k,0))\right)$.
See figure~~\ref{cyl_with_rect}.

\begin{cor}    \label{cyl_with_rect_cor}
The set of parameters $\het\in\R/2\pi\Z$ such that the noncompact
polygon $\tP_{\het}$ is conservative is a dense $G_{\de}$.
\begin{proof}
In view of Theorem~~\ref{mor_dens_gdel_thm}, it suffices to check
that $(C_0,\ka)$ is an amenable pair. The surface $P_{\het}$ is
rational iff $\frac{\het}{\pi}=\frac{m}{n}$, with $m$ and $n$ relatively prime. We then have
$N(P_{\het})=n$ if $n$ is even, and $N(P_{\het})=2n$ if $n$ is odd. Thus,
$N(P_{\het})$ is even for any $\pi$-rational $\het$.
The claim now follows from the density of rational numbers.
\end{proof}
\end{cor}

We expect that among $\Z$-periodic polygonal surfaces, conservativeness is generic
in a stronger sense than that ensured by Theorem~~\ref{dens_gdel_recur_thm} and Theorem~~\ref{mor_dens_gdel_thm}.
For instance, let $\{\tP_{\het}:\het\in\R/2\pi\Z\}$ be the family in Corollary~~\ref{cyl_with_rect_cor}.
We conjecture that the noncompact surfaces $\tP_{\het}$ are conservative for all $\het$.
However, as the following example shows, not every $\Z$-periodic polygonal surface with a boundary is conservative.

\begin{figure}[htbp]
\begin{center}
\input{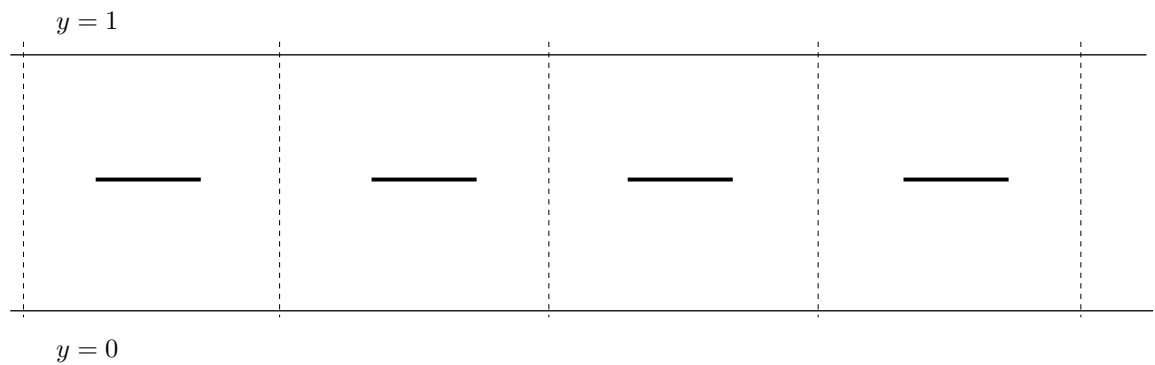}
\caption{Infinite band with a $\Z$-periodic collection of horizontal
barriers.} \label{nonconserv}
\end{center}
\end{figure}

\begin{exa}    \label{nonconserv_exa}
{\em Let $0<l<1$. Let $\tP$ be the standard band with a $\Z$-periodic
collection of horizontal barriers of length $l$. See figure~~\ref{nonconserv} where
the barriers are in the middle of the band. This surface belongs to the family of surfaces
considered in Corollary~~\ref{band_with_bar_cor}. Also, $\tP$ is a degeneration of the surface
$\tP_0$ in the family $\{\tP_{\het}:\het\in\R/2\pi\Z\}$ of Corollary~~\ref{cyl_with_rect_cor}.
The quotient $P=\tP/\Z$ is the standard cylinder with a horizontal barrier. It is a rational polygonal surface;
$N(P)=1$. Thus, Theorem~~\ref{per_rat_recurr_thm} does not apply.

The horizontal component of a tangent vector does not change under the geodesic flow
of $\tP$. Hence, the geodesic flow is transient. The directional flows $(U\tP_{\het},\tG_{\het}^t,\tmu_{\het})$ are
also transient for $\het\ne 0$, while the flow $\tG_0^t$ is periodic.

The fact that our barriers are positioned at the height $1/2$ plays no role
in this case. Any $\Z$-periodic  collection of horizontal barriers in $B_0$
leads to a transitive polygonal surface. We leave further generalizations to the reader.

}
\end{exa}

\end{document}